\documentstyle[12pt,leqno]{article}
\font\erls=eufm10 at 11pt

\newcommand{\Ls}{\mbox{\erls L}}
\font\erl=eufm10 at 12pt
\newcommand{\M}{\mbox{\erl M}}
\newcommand{\LL}{\mbox{\erl L}}
\font\erll=eufm10 at 15pt
\newcommand{\LLL}{\mbox{\erll L}}
\font\ers=eufm10 
\font\ers=eusm10 at 12pt

\newcommand{\HH}{\mbox{\ers H}}
\font\erss=eusm9 at 9pt

\newcommand{\HHs}{\mbox{\erss H}}

\font\ms=msam8 at 11pt
\newcommand{\eql}{\mbox{\thinspace\thinspace{\ms 5}\thinspace\thinspace}}
\newcommand{\eqg}{\mbox{\thinspace\thinspace{\ms =}\thinspace\thinspace}}
\font\mss=msam8 at 9pt
\newcommand{\eqls}{\mbox{{\mss 5}}}

\font\msb=msbm8 at 12pt
\newcommand{\R}{\mbox{\msb R}}
\newcommand{\C}{\mbox{\msb C}}

\newcommand{\Rs}{\mbox{\msb R}}

\begin{document}
\vspace*{25mm}
\begin{center}
{\large \bf Integral Transforms with{\boldmath$H$}-Function Kernels on 
$\LLL_{\nu,r}$-Spaces}\\[5mm]
Hans-J\"urgen {\sc Glaeske}\\
{\small \it Department of Mathematics\\Friedrich Schiller University$,$ 
D-07740 Jena$,$ Germany}\\[5mm]
Anatoly A. {\sc Kilbas}\\
{\small \it Department of Mathematics and Mechanics\\Belarusian State 
University$,$ Minsk 220050$,$ Belarus}\\[5mm]
Megumi {\sc Saigo}\\
{\small \it Department of Applied Mathematics\\Fukuoka University$,$ 
Fukuoka 814-0180$,$ Japan}\\[5mm]
Sergei A. {\sc Shlapakov}\\
{\small \it Department of Mathematics\\Vitebsk Pedagogical University$,$ 
Vitebsk 210036$,$ Belarus}\\[8mm]
\end{center}
\begin{abstract}
Integral transforms
$$(\mbox{\boldmath$H$}f)(x)=\int^\infty_0H^{m,n}_{\thinspace p,q}
\left[xt\left|\begin{array}{c}(a_i,\alpha_i)_{1,p}\\[1mm](b_j,\beta_j)_{1,q}
\end{array}\right.\right]f(t)dt$$
involving Fox's $H$-functions as kernels are studied in the spaces 
$\Ls_{\nu,r}$ of functions $f$ such that
$$\int^\infty_0|t^\nu f(t)|^r\frac{dt}t<\infty\quad(1\ \eqls\ r<\infty,
\ \nu\in\Rs).$$
Mapping properties such as the boundedness, the representation and the range 
of the transforms \boldmath$H$ are given.
\end{abstract}
\begin{flushleft}
{\bf 1. \ Introduction}
\end{flushleft}
\setcounter{section}{1}
\setcounter{equation}{0}
\par
This paper deals with the integral transforms of the form
\begin{eqnarray}
(\mbox{\boldmath$H$}f)(x)=\int^\infty_0H^{m,n}_{\thinspace p,q}
\left[xt\left|\begin{array}{c}(a_i,\alpha_i)_{1,p}\\[2mm](b_j,\beta_j)_{1,q}
\end{array}\right.\right]f(t)dt,
\end{eqnarray}
where $ H^{m,n}_{\thinspace p,q}\left[z\left|\begin{array}{c}
(a_i,\alpha_i)_{1,p}\\[2mm](b_j,\beta_j)_{1,q}\end{array}\right.\right]$ is 
the Fox $H$-function. This function of general hypergeometric type was 
introduced by Fox [8]. For integers $m,n,p,q$ such that $0\eql m\eql q,$ 
$0\eql n\eql p$, $a_i,b_j\in\C$ with \C\ of the field of complex numbers, and 
$\alpha_i,\beta_j\in\R_+=(0,\infty)\ (1\eql i\eql p,$ $1\eql j\eql q)$, it is 
defined by
\begin{eqnarray}
H^{m,n}_{\thinspace p,q}\left[z\left|\begin{array}{c}(a_i,\alpha_i)_{1,p}\\[2mm]
(b_j,\beta_j)_{1,q}\end{array}\right.\right]&\hspace{-2.5mm}=&\hspace{-2.5mm}
H^{m,n}_{\thinspace p,q}\left[z\left|\begin{array}{ll}(a_1,\alpha_1),
&\cdots,(a_p,\alpha_p)\\[2mm](b_1,\beta_1),&\cdots,(b_q,\beta_q)\end{array}
\right.\right]\\[2mm]
       &\hspace{-2.5mm}=&\hspace{-2.5mm}{\frac1{2\pi i}}\int_{\Ls}
\HH^{m,n}_{\thinspace p,q}\left[\left.\begin{array}{c}(a_i,\alpha_i)_{1,p}
\\[2mm](b_j,\beta_j)_{1,q}\end{array}\right|s\right]z^{-s}ds,\nonumber
\end{eqnarray}
where
\begin{eqnarray}
\HH^{m,n}_{\thinspace p,q}\left[\left.\begin{array}{c}(a_i,\alpha_i)_{1,p}
\\[2mm](b_j,\beta_j)_{1,q}\end{array}\right|s\right]
={\frac{\displaystyle{\prod^m_{j=1}\Gamma(b_j+\beta_js)\prod^n_{i=1}
\Gamma(1-a_i-\alpha_is)}}{\displaystyle{\prod^p_{i=n+1}\Gamma(a_i+\alpha_is)
\prod^q_{j=m+1}\Gamma(1-b_j-\beta_js)}}},
\end{eqnarray}
the contour $\LL$ is specially chosen and an empty product, if it occurs, is 
taken to be one. The theory of this function may be found in [2], [29, Chapter 
1], [36, \S8.3] and [52, Chapter 2]. We abbreviate the Fox $H$-function (1.2) 
and the function in (1.3) to $H^{m,n}_{\thinspace p,q}(x)$, 
$\HH^{m,n}_{\thinspace p,q}(s)$ or $H(x)$, $\HH(s)$ when no confusion occurs.
\par
Most of the known integral transforms can be put into the form (1.1). When 
$\alpha_1=\cdots=\alpha_p=\beta_1=\cdots=\beta_q=1$, then (1.2) is the Meijer 
$G$-function [7, Chapter 5.3] and (1.1) is reduced to the so-called integral 
transforms with $G$-function kernels or {\boldmath$G$}-transforms. The 
classical Laplace and Hankel transforms, the Riemann-Liouville fractional 
integrals, the even and odd Hilbert transforms, the integral transforms with 
the Gauss hypergeometric function, etc. belong to these 
{\boldmath$G$}-transforms, for whose theory and historical notices see 
[44, \S\S36, 39]. There are other transforms which can not be reduced to 
{\boldmath$G$}-transforms but can be put into the transform {\boldmath$H$} 
given in (1.1). Such kinds of transforms are the modified Laplace and Hankel 
transforms [50], [52], [40], [44, \S\S18, 23, 39], the Erd\'elyi-Kober type 
fractional integration operators [26], [6], [50], [44, \S18], the transforms 
with the Gauss hypergeometric function as kernels [35], [30], [41], [42], 
[44, \S\S23, 39], the Bessel-type integral transforms [27], [38], [21], [22], 
etc.
\par
The integral transforms (1.1) with $H$-function kernels or 
{\boldmath$H$}-transforms were first considered by Fox [8] while investigating 
$G$- and $H$-functions as symmetrical Fourier kernels. This paper together 
with the ones [17], [45], [9], [10], [49], [14], [4], [28] and [34] were 
devoted to find the inversion formulae for the {\boldmath$H$}-transforms (1.1) 
is the spaces $L_1(0,\infty)$ and $L_2(0,\infty)$. Some properties of 
{\boldmath$H$}-transforms such as their Mellin transform, the relation of 
fractional integration by parts, compositional formulae, etc. were considered 
in [11], [12], [13], [51], [46] and [15]. In [47], [5] and [1] the operators 
{\boldmath$H$} in (1.1) were represented as the compositions of the 
Erd\'elyi-Kober type operators and the integral operators of the form (1.1) 
with the $H$-function of the less order. Factorization properties of (1.1) in 
special functional spaces $L^\Phi_2$ were investigated in [53]. The properties 
of generalized fractional integration operators, being some modifications of 
the operators (1.1), with $H^{m,0}_{\thinspace m,m}$-function as the kernel 
were investigated on the space $L_p(0,\infty)$ in [24], [16], [25] and on 
McBride spaces $F_{p,\mu}$ and $F_{p,\mu}'$ (see [31] and [44, \S8]) in [37], 
[43]. We also note that some facts about multidimensional transforms 
{\boldmath$H$} of the form (1.1) were given in [3, \S4.4].
\par
The papers [32], [33] were devoted to the range and the invertibility of the 
operators (1.1) in the special cases when $\HH(s)=\Gamma(b+s/m),
\Gamma(1-a-s/m)$ and $\Gamma(b+s/m)/\Gamma(1-a-s/m)$ with $m>0$, in the spaces 
$\LL_{\nu,r}$ and its subspaces $F_{p,\mu}$. The former space with 
$\nu\in\R=(-\infty,\infty)$ and $1\eql r<\infty$ is defined as the space of 
those Lebesgue measurable complex functions $f$ for which
\begin{eqnarray}
\|f\|_{\nu,r}=\int^\infty_0\left|t^\nu f(t)\right|^r\frac{dt}t<\infty
\end{eqnarray}
(see [40]). Our previous papers [18], [23] were dealt with the study of the 
mapping properties such as the boundedness, the representation and the range 
of the general transforms {\boldmath$H$} defined by (1.1) in the space 
$\LL_{\nu,2}$. The results in [18], [23] were extended to the space 
$\LL_{\nu,r}$ with any $1\eql r<\infty$ in [19], [20], [48]. It is proved that 
the obtained results are different in nine cases:
\begin{eqnarray*}
&\begin{array}{lll}
1) \ a^*=\Delta={\rm Re}(\mu)=0;&2) \ a^*=\Delta=0,\ {\rm Re}(\mu)<0;
&3) \ a^*=0,\ \Delta>0;\\[2mm]
4) \ a^*=0,\ \Delta<0;&5) \ a_1^*>0,\ a_2^*>0;&6) \ a_1^*>0,\ a_2^*=0;\\[2mm]
7) \ a_1^*=0,\ a_2^*>0;&8) \ a^*>0,\ a_1^*>0,\ a_2^*<0;
&9) \ a^*>0,\ a_1^*<0,\ a_2^*>0.
\end{array}&
\end{eqnarray*}
Here
\begin{eqnarray}
a^*&\hspace{-2.5mm}=&\hspace{-2.5mm}\sum^n_{i=1}\alpha_i-\sum^p_{i=n+1}
\alpha_i+\sum^m_{j=1}\beta_j-\sum^q_{j=m+1}\beta_j,\\[2mm]
\Delta&\hspace{-2.5mm}=&\hspace{-2.5mm}\sum^q_{j=1}\beta_j-\sum^p_{i=1}
\alpha_i,\\[2mm]
\mu&\hspace{-2.5mm}=&\hspace{-2.5mm}\sum^q_{j=1}b_j-\sum^p_{i=1}a_i
+\frac{p-q}2,\\[2mm]
a_1^*&\hspace{-2.5mm}=&\hspace{-2.5mm}\sum^m_{j=1}\beta_j-\sum^p_{i=n+1}
\alpha_i,\qquad a_2^*=\sum^n_{i=1}\alpha_i-\sum^q_{j=m+1}\beta_j.
\end{eqnarray}
We note that
\begin{eqnarray}
&&a^*=a_1^*+a^*_2,\qquad\Delta=a_1^*-a_2^*.
\end{eqnarray}
The results in [19], [20] and [48] being the extensions of those by Rooney 
[40] from {\boldmath$G$}-transforms to {\boldmath$H$}-transforms were proved 
under the additional condition $\delta=1$, where
\begin{eqnarray}
&&\delta=\prod^p_{i=1}\alpha_i^{-\alpha_i}\prod^q_{j=1}\beta_j^{\beta_j}.
\end{eqnarray}
\par
The present paper is devoted to extend the results in [19], [20], [48] from 
$\delta=1$ to any $\delta>0$. Section 2 contains preliminary information from 
[40] concerning the mapping properties on $\LL_{\nu,r}$ of the Mellin 
transform, Erd\'elyi-Kober type fractional integral operators and generalized 
Hankel and Laplace transforms. In Section 3 we summarize the results from [18] 
on the asymptotic properties of the function $\HH(s)$ defined in (1.3) and its 
derivative and on $\LL_{\nu,2}$-theory of the transform {\boldmath$H$}. 
Sections 4 and 5 deal with the one-to-one boundedness, representation and 
range of {\boldmath$H$}-transform (1.1) on the space 
$\LL_{\nu,r}\ (1\eql r<\infty)$ in the cases when $a^*=\Delta=0$ and 
$a^*=0,\ \Delta\ne0$, respectively. Sections 6 and 7 are devoted to consider 
the cases $a_1^*\eqg0,\ a_2^*\eqg0$ and $a^*>0,\ a_1^*<0$ or $a_2^*<0$, 
respectively.
\\\\
\begin{flushleft}
{\bf 2. \ Some Auxiliary Results}
\end{flushleft}
\setcounter{section}{2}
\setcounter{equation}{0}
\par
In this section we collect a variety of facts concerning multipliers for the 
Mellin transform and well known integral operators (see [39], [40]) which we 
need in next sections. For $f\in \LL_{\nu,r}$ with $1\eql r\eql2$, the Mellin 
transform of $f$ is defined [40] by
\begin{eqnarray}
&&\left(\M f\right)(\nu+it)=\int^\infty_{-\infty}e^{(\nu+it)\tau}f(e^\tau)d\tau
\end{eqnarray}
for $\nu,t\in\R$. We also write $\left(\M f\right)(s)$ for ${\rm Re}(s)=\nu$ 
as $\left(\M f\right)(\nu+it)$. In particular, if $f\in \LL_{\nu,r}\bigcap
\LL_{\nu,1}$, then $\left(\M f\right)(s)$ is given by the usual expression
\begin{eqnarray}
&&\left(\M f\right)(s)=\int^\infty_0f(t)t^{s-1}dt.
\end{eqnarray}
\par
First we give the definition of the set ${\cal A}$ and formulate the 
multiplier theorem for the Mellin transform \M.
\\\par
{\bf Definition 1.} \ [39, Definition 3.1] \ We say that a function $m$ 
belongs to ${\cal A}$ if there are extended real numbers $\alpha(m)$ and 
$\beta(m)$ with $\alpha(m)<\beta(m)$ such that:
\begin{enumerate}
\item[{\bf a)\ }] \ $m(s)$ is analytic in the strip $\alpha(m)<
{\rm Re}(s)<\beta(m)$;
\item[{\bf b)\ }] \ $m(s)$ is bounded in every closed substrip 
$\sigma_1\eql{\rm Re}(s)\eql\sigma_2$, where $\alpha(m)<\sigma_1\eql
\sigma_2<\beta(m)$;
\item[{\bf c)\ }] \ $|m'(\sigma+it)|=O\left(|t|^{-1}\right)$ as 
$|t|\to\infty$ for $\alpha(m)<\sigma<\beta(m)$.
\\
\end{enumerate}
\par
For two Banach spaces $X$ and $Y$ we use the notation $[X,Y]$ denoting the 
collection of bounded linear operators from $X$ to $Y$, and $[X,X]$ is 
abbreviated to $[X]$.
\\\par
{\bf Theorem 1.} \ [39, Theorem 1] \ {\sl Let $m\in{\cal A},
\alpha(m)<\nu<\beta(m)$ and $1<r<\infty$. Then there is a transform 
$T_m\in[\LL_{\nu,r}]$ such that for every $f\in \LL_{\nu,r}$ with $1<r\eql2,$ 
the relation
\begin{eqnarray}
&&\left(\M T_mf\right)(s)=m(s)\left(\M f\right)(s)\quad({\rm Re}(s)=\nu)
\end{eqnarray}
holds valid. For $\alpha(m)<\nu<\beta(m)$ and $1<r\eql2$ the transform $T_m$ 
is one-to-one on $\LL_{\nu,r}$ except when $m=0$. If $1/m\in{\cal A},$ then 
for $\max[\alpha(m),\alpha(1/m)]<\nu<\min[\beta(m),\beta(1/m)]$ and for 
$1<r<\infty,$ $T_m$ maps $\LL_{\nu,r}$ one-to-one onto itself$,$ and there 
holds}
\begin{eqnarray*}
&&T_m^{-1}=T_{1/m}.
\end{eqnarray*}
\\\par
In the discussion of this article we use the special integral operators as 
follows: The {\it Erd\'elyi-Kober type fractional integrals} (see [44, \S18.1]) 
for $\alpha,\eta\in\C$ with ${\rm Re}(\alpha)>0$ and $\sigma,x\in\R_+$:
\begin{eqnarray}
&&\left(I^\alpha_{0+;\sigma,\eta}f\right)(x)=\frac{\sigma x^{-\sigma(\alpha+
\eta)}}{\Gamma(\alpha)}\int^x_0(x^\sigma-t^\sigma)^{\alpha-1}t^{\sigma\eta+
\sigma-1}f(t)dt;\\[2mm]
&&\left(I^\alpha_{-;\sigma,\eta}f\right)(x)=\frac{\sigma x^{\sigma\eta}}
{\Gamma(\alpha)}\int^\infty_x(t^\sigma-x^\sigma)^{\alpha-1}
t^{\sigma(1-\alpha-\eta)-1}f(t)dt;
\end{eqnarray}
the {\it modified Hankel transform} for $\kappa\in\R\backslash\{0\},
{\rm Re}(\eta)>-1$ and $x\in\R_+$:
\begin{eqnarray}
&&\left(H_{\kappa,\eta}f\right)(x)=\int^\infty_0(xt)^{1/\kappa-1/2}
J_\eta\left(|\kappa|(xt)^{1/\kappa}\right)f(t)dt;
\end{eqnarray}
and the {\it modified Laplace transform} for $\kappa\in\R\backslash\{0\},
\alpha\in\C$ and $x\in\R_+$:
\begin{eqnarray}
&&\left(L_{\kappa,\alpha}f\right)(x)=\int^\infty_0(xt)^{-\alpha}
e^{-|\kappa|(xt)^{1/\kappa}}f(t)dt.
\end{eqnarray}
All these transforms are defined for continuous functions $f$ with compact 
support on $\R_+$ for the range of parameters indicated.
\par
For $\nu \in \R_+$ we denote by $\LL_{\nu,\infty}$ the collection of functions 
$f$, measurable on $\R_{+}$, such that
$$\|f\|_{\nu,\infty}=\mathop{\rm ess\thinspace sup}_{x>0}|x^\nu f(x)|<\infty.$$
\par
The boundedness properties of the transforms (2.4) - (2.7) and their Mellin 
transforms are given by the following statement.
\\\par
{\bf Theorem 2.} \ [40, Theorem 5.1] \ {\bf (a)} \ {\sl If 
$1\eql r\eql\infty,{\rm Re}(\alpha)>0$ and $\nu<\sigma(1+{\rm Re}(\eta)),$ 
then for all $s\eqg r$ such that $1/s>1/r-{\rm Re}(\alpha),$ the operator 
$I^\alpha_{0+;\sigma,\eta}$ belongs to $[\LL_{\nu,r},\LL_{\nu,s}]$ and is a 
one-to-one transform from $\LL_{\nu,r}$ onto $\LL_{\nu,s}$. For $1\eql r\eql2$ 
and $f\in \LL_{\nu,r}$
\begin{eqnarray}
&&\left(\M I^\alpha_{0+;\sigma,\eta}f\right)(s)=\frac{\Gamma(1+\eta-s/\sigma)}
{\Gamma(1+\eta+\alpha-s/\sigma)}\left(\M f\right)(s)\quad({\rm Re}(s)=\nu).
\end{eqnarray}
{\bf (b)} \ If $1\eql r\eql\infty,{\rm Re}(\alpha)>0$ and 
$\nu>-\sigma{\rm Re}(\eta),$ then for all $s\eqg r$ such that 
$1/s>1/r-{\rm Re}(\alpha),$ the operator $I^\alpha_{-;\sigma,\eta}$ belongs to 
$[\LL_{\nu,r},\LL_{\nu,s}]$ and is a one-to-one transform from $\LL_{\nu,r}$ 
onto $\LL_{\nu,s}$. For $1\eql r\eql2$ and $f\in \LL_{\nu,r}$
\begin{eqnarray}
&&\left(\M I^\alpha_{-;\sigma,\eta}f\right)(s)=\frac{\Gamma(\eta+s/\sigma)}
{\Gamma(\eta+\alpha+s/\sigma)}\left(\M f\right)(s)\quad({\rm Re}(s)=\nu).
\end{eqnarray}
{\bf (c)} \ If $1<r<\infty$ and $\gamma(r)\eql\kappa(\nu-1/2)+1/2<
{\rm Re}(\eta)+3/2,$ where
\begin{eqnarray}
&&\gamma(r)=\max\left[\frac1r,\frac1{r'}\right]\quad\mbox{for}\quad
\frac1r+\frac1{r'}=1,
\end{eqnarray}
then for all $s\eqg r$ such that $s'\eqg(\kappa(\nu-1/2)+1/2)^{-1}$ and 
$1/s+1/s'=1,$ the operator $H_{\kappa,\eta}$ belongs to 
$[\LL_{\nu,r},\LL_{1-\nu,s}]$ and is a one-to-one transform from $\LL_{\nu,r}$ 
onto $\LL_{1-\nu,s}$ If $1<r\eql2$ and $f\in \LL_{\nu,r},$ then
\begin{eqnarray}
&&\left(\M H_{\kappa,\eta}f\right)(s)=\left(\frac2{|\kappa|}
\right)^{\kappa(s-1/2)}\frac{\Gamma([\eta+\kappa(s-1/2)+1]/2)}
{\Gamma([\eta-\kappa(s-1/2)+1]/2)}\left(\M f\right)(1-s)\\[2mm]
&&\hspace{105mm}({\rm Re}(s)=1-\nu).\nonumber
\end{eqnarray}
{\bf (d)} \ If $1\eql r\eql s\eql\infty,$ and if $\nu<1-{\rm Re}(\alpha)$ for 
$\kappa>0$ and $\nu>1-{\rm Re}(\alpha)$ for $\kappa<0,$ then the operator 
$L_{\kappa,\alpha}$ belongs to $[\LL_{\nu,r},\LL_{1-\nu,s}]$ and is a 
one-to-one transform from $\LL_{\nu,r}$ onto $\LL_{1-\nu,s}$. If $1\eql r\eql2$ 
and $f\in \LL_{\nu,r},$ then}
\begin{eqnarray}
&&\left(\M L_{\kappa,\alpha}f\right)(s)=\Gamma(\kappa[s-\alpha])
|\kappa|^{1-\kappa(s-\alpha)}\left(\M f\right)(1-s)\quad({\rm Re}(s)=1-\nu).
\end{eqnarray}
\\\par
For further investigation we also need certain elementary operators. For a 
function $f$ being defined almost everywhere on $\R_+$ we denote the operators 
$M_\zeta $, $W_\delta $ and $R$ as follows:
\begin{eqnarray}
&&\left(M_\zeta f\right)(s)=x^\zeta f(x)\quad\mbox{for }\zeta\in\C;\\[2mm]
&&\left(W_\delta f\right)(s)=f\left(\frac x\delta\right)\quad
\mbox{for }\delta\in\R_+;\\[2mm]
&&\left(Rf\right)(x)=\frac1x f\left(\frac 1x\right).
\end{eqnarray}
\par
These operators have the following properties for $\nu\in\R$ and 
$1\eql r<\infty$ (see [40]):
\par
{\bf (P1)} \ $M_\zeta $ is an isometric isomorphism of $\LL_{\nu,r}$ onto 
$\LL_{\nu-{\rm Re}(\zeta),r}$, and if $f\in \LL_{\nu,r}$ $(1\eql r\eql2)$, 
then
\begin{eqnarray}
&&\left(\M M_\zeta f\right)(s)=\left(\M f\right)(s+\zeta)\quad
({\rm Re}(s)=\nu-{\rm Re}(\zeta));
\end{eqnarray}
\par
{\bf (P2)} \ $W_\delta $ is an isometric isomorphism of $\LL_{\nu,r}$, and if 
$f\in \LL_{\nu,r}$ $(1\eql r\eql2)$, then
\begin{eqnarray}
&&\left(\M W_\delta f\right)(s)=\delta^s\left(\M f\right)(s)\quad
({\rm Re}(s)=\nu);
\end{eqnarray}
\par
{\bf (P3)} \ $R$ is an isometric isomorphism of $\LL_{\nu,r}$ onto 
$\LL_{1-\nu,r}$, and if $f\in \LL_{\nu,r}$ $(1\eql r\eql2)$, then
\begin{eqnarray}
&&\left(\M R f\right)(s)=\left(\M f\right)(1-s)\quad({\rm Re}(s)=1-\nu).
\end{eqnarray}
\\\\
\newpage
\begin{flushleft}
{\bf 3. \ Asymptotic Properties of $\HH^{m,n}_{p,q}(s)$ and 
$\LL_{\nu,2}$-Theory of {\boldmath$H$}-Transforms}
\end{flushleft}
\setcounter{section}{3}
\setcounter{equation}{0}
\par
Let $a^*,\Delta,\mu,a_1^*,a_2^*$ and $\delta$ be given by (1.5) - (1.8) and 
(1.10), respectively. For integers $m,n,p,q$ such that 
$0\eql m\eql q,0\eql n\eql p$ and for $a_1,\cdots,a_p,b_1,\cdots,b_q\in\C$ and 
for $\alpha_1,\cdots,\alpha_p,\beta_1,\cdots,\beta_q\in\R_+$ we define $\alpha$ 
and $\beta$ by
\begin{eqnarray}
&\alpha=\left\{\begin{array}{ll}\displaystyle\max\left[-\frac{{\rm Re}(b_1)}
{\beta_1},\cdots,-\frac{{\rm Re}(b_m)}{\beta_m}\right]&\mbox{if}\ m>0,\\[4mm]
-\infty&\mbox{if}\ m=0;\end{array}\right.&\\[2mm]
&\beta=\left\{\begin{array}{ll}\displaystyle\min\left[\frac{1-{\rm Re}(a_1)}
{\alpha_1},\cdots,\frac{1-{\rm Re}(a_n)}{\alpha_n}\right]&\mbox{if}\ n>0,\\[4mm]
\infty&\mbox{if}\ n=0.\end{array}\right.&
\end{eqnarray}
\\\par
{\bf Lemma 1.} \ [18, Lemma 1] \ {\sl If $\sigma,t\in\R,$ then the estimate
\begin{eqnarray}
\left|\HH^{m,n}_{\thinspace p,q}(\sigma+it)\right|&\hspace{-2.5mm}\sim
&\hspace{-2.5mm}(2\pi)^{c^*}\delta^\sigma\prod^p_{i=1}\alpha_i^{1/2-{\rm Re}
(a_i)}\prod^q_{j=1}\beta_j^{{\rm Re}(b_j)-1/2}\\[2mm]
&&\times\ |t|^{\Delta\sigma+{\rm Re}(\mu)}\exp\left\{-\frac{\pi|t|a^*}2
-\frac{\pi{\rm Im}(\xi){\rm sign}(t)}2\right\}\nonumber
\end{eqnarray}
with
\begin{eqnarray}
&&c^*=m+n-\frac{p+q}2,\qquad\xi=\sum^n_{i=1}a_i-\sum^p_{i=n+1}a_i+
\sum^m_{j=1}b_j-\sum^q_{i=m+1}b_j
\end{eqnarray}
holds as $|t|\to\infty$ uniformly in $\sigma$ for $\sigma$ in any bounded 
interval in \R. Further},
\begin{eqnarray}
\left\{\HH^{m,n}_{\thinspace p,q}(\sigma+it)\right\}'&\hspace{-2.5mm}=
&\hspace{-2.5mm}\HH^{m,n}_{\thinspace p,q}(\sigma+it)\Biggl[\log\delta+a^*_1
\log(it)\\[2mm]
&&-a^*_2\log(-it)+\frac{\mu+\Delta\sigma}{it}+O\left(\frac1{t^2}\right)\Biggr]
\quad(|t|\to\infty).\nonumber
\end{eqnarray}
\\\par
To present $\LL_{\nu,2}$-theory for {\boldmath$H$}-transform we have to define 
a certain set of real numbers.
\\\par
{\bf Definition 2.} \ For the function $\HH(s)$ given in (1.3) we call the 
{\it exceptional set} ${\cal E}_{\HHs}$ of {\HH} the set of real numbers $\nu$ 
such that $\alpha<1-\nu<\beta$ and $\HH(s)$ has a zero on the line 
${\rm Re}(s)=1-\nu$.
\\\par
{\bf Theorem 3.} \ [18, Theorem 3] \ {\sl Suppose that
\begin{enumerate}
\item[{{\bf (a)}\ }]$\alpha<1-\nu<\beta$
\end{enumerate}
and that either of the conditions
\begin{enumerate}
\item[{{\bf (b)}\ }]$a^*>0;$
\item[{{\bf (c)}\ }]$a^*=0,\Delta(1-\nu)+{\rm Re}(\mu)\eql0.$
\end{enumerate}
holds$,$ then we have:
\par
{\bf (i)} \ There is a one-to-one transform $\mbox{\boldmath$H$}
\in[\LL_{\nu,2},\LL_{1-\nu,2}]$ so that the relation
\begin{eqnarray}
&&\left(\M\mbox{\boldmath$H$}f\right)(s)=\HH^{m,n}_{\thinspace p,q}
\left[\left.\begin{array}{c}(a_p,\alpha_p)\\[2mm](b_q,\beta_q)\end{array}
\right|s\right]\left(\M f\right)(1-s)
\end{eqnarray}
holds for $f\in\LL_{\nu,r}$ and ${\rm Re}(s)=1-\nu$. If 
$a^*=0,\ \Delta(1-\nu)+{\rm Re}(\mu)=0$ and $\nu\notin{\cal E}_{\HHs},$ then 
the operator {\boldmath$H$} maps $\LL_{\nu,2}$ onto $\LL_{1-\nu,2}$.
\par
{\bf (ii)} \ For $f,g\in \LL_{\nu,2}$ the relation 
\begin{eqnarray}
&&\int^\infty_0f(x)\left(\mbox{\boldmath$H$}g\right)(x)dx=\int^\infty_0g(x)
\left(\mbox{\boldmath$H$}f\right)(x)dx
\end{eqnarray}
holds.
\par
{\bf (iii)} \ Let $f\in \LL_{\nu,2},\lambda\in\C$ and $h>0$. If 
${\rm Re}(\lambda)>(1-\nu)h-1,$ then $\mbox{\boldmath$H$}f$ is given by
\begin{eqnarray}
&&\left(\mbox{\boldmath$H$}f\right)(x)=hx^{1-(\lambda+1)/h}\frac d{dx}
\ x^{(\lambda+1)/h}\\[2mm]
&&\hspace{10mm}\times\ \int^\infty_0H^{m,n+1}_{\thinspace p+1,q+1}
\left[xt\left|\begin{array}{c}(-\lambda,h),(a_1,\alpha_1),\cdots,(a_p,\alpha_p)
\\[2mm](b_1,\beta_1)\cdots,(b_q,\beta_q),(-\lambda-1,h)\end{array}\right.
\right]f(t)dt.\nonumber
\end{eqnarray}
If ${\rm Re}(\lambda)<(1-\nu)h-1,$ then
\begin{eqnarray}
&&\left(\mbox{\boldmath$H$}f\right)(x)=-hx^{1-(\lambda+1)/h}
\frac d{dx}\ x^{(\lambda+1)/h}\\[2mm]
&&\hspace{10mm}\times\ \int^\infty_0H^{m+1,n}_{\thinspace p+1,q+1}
\left[xt\left|\begin{array}{c}(a_1,\alpha_1),\cdots,(a_p,\alpha_p),(-\lambda,h)
\\[2mm](-\lambda-1,h),(b_1,\beta_1)\cdots,(b_q,\beta_q)\end{array}\right.
\right]f(t)dt.\nonumber
\end{eqnarray}
\par
{\bf (iv)} \ {\boldmath$H$} is independent on $\nu$ in the sense that if 
$\nu_1$ and $\nu_2$ satisfy {\rm (a),} and {\rm (b)} or {\rm (c),} and if the 
transforms $\mbox{\boldmath$H$}_1$ and $\mbox{\boldmath$H$}_2$ are given by 
$(3.6),$ then $\mbox{\boldmath$H$}_1f=\mbox{\boldmath$H$}_2f$ for} 
$f\in \LL_{\nu_1,2}\bigcap \LL_{\nu_2,2}$.
\\\par
{\bf Corollary.} \ {\sl Let $\alpha<\beta$ and one of the following conditions 
holds$:$
\begin{enumerate}
\item[{{\bf (b)}\ }]$a^*>0;$
\item[{{\bf (e)}\ }]$a^*=0,\Delta>0$ and $\displaystyle\alpha<
-\frac{{\rm Re}(\mu)}\Delta;$
\item[{{\bf (f)}\ }]$a^*=0,\Delta<0$ and $\displaystyle\beta>
-\frac{{\rm Re}(\mu)}\Delta;$
\item[{{\bf (g)}\ }]$a^*=0,\Delta=0$ and ${\rm Re}(\mu)\eql0.$
\end{enumerate}
Then the transform {\boldmath$H$} can be defined on $\LL_{\nu,2}$ with} 
$1-\beta<\nu<1-\alpha$.
\\\par
{\bf Theorem 4.} \ [18, Theorem 4] \ {\sl Let $\alpha<1-\nu<\beta$ and either 
of the following conditions holds$:$
\begin{enumerate}
\item[{{\bf (b)}\ }]$a^*>0;$
\item[{{\bf (d)}\ }]$a^*=0,\Delta(1-\nu)+{\rm Re}(\mu)<-1$.
\end{enumerate}
Then for $x\in\R_+,\left(\mbox{\boldmath$H$}f\right)(x)$ is given by $(1.1)$ 
for} $f\in \LL_{\nu,2}.$
\\\par
{\bf Corollary.} \ {\sl Let $\alpha<\beta$ and one of the following conditions 
holds$:$
\begin{enumerate}
\item[{{\bf (b)}\ }]$a^*>0;$
\item[{{\bf (h)}\ }]$a^*=0,\Delta>0$ and $\displaystyle\alpha<
-\frac{1+{\rm Re}(\mu)}\Delta;$
\item[{{\bf (i)}\ }]$a^*=0,\Delta<0$ and $\displaystyle\beta>
-\frac{1+{\rm Re}(\mu)}\Delta;$
\item[{{\bf (j)}\ }]$a^*=0,\Delta=0$ and ${\rm Re}(\mu)<-1.$
\end{enumerate}
Then the transform {\boldmath$H$} can be defined by $(1.1)$ on $\LL_{\nu,2}$ 
with} $1-\beta<\nu<1-\alpha$.
\\\\
\begin{flushleft}
{\bf 4. \ $\LL_{\nu,r}$-Theory of the Transform {\boldmath$H$} ($a^*=\Delta=0$)}
\end{flushleft}
\setcounter{section}{4}
\setcounter{equation}{0}
\par
In this section, basing on the existence of the transform {\boldmath${H}$} on 
the space $\LL_{\nu,2}$ which is guaranteed in Theorem 3 for some $\nu\in\R$ 
and $a^*=\Delta=0$, we prove that such a transform can be extended to 
$\LL_{\nu,r}$ for $1<r<\infty$ such that 
{\boldmath${H}$}$\in[\LL_{\nu,r},\LL_{1-\nu,s}]$ for a certain range of the 
value $s$. We also characterize the range of {\boldmath${H}$} on $\LL_{\nu,r}$ 
in terms of the Erd\'elyi-Kober type fractional integral operators 
$I^\alpha_{0+;\sigma,\eta}$ and $I^\alpha_{-;\sigma,\eta}$ given in (2.4) and 
(2.5) except for its isolated values $\nu\in{\cal E}_{\HHs}$. The results will 
be different in the cases ${\rm Re}(\mu)=0$ and ${\rm Re}(\mu)\ne0$, where 
$\mu$ is defined by (1.7). First we consider the former case.
\\\par
{\bf Theorem 5.} \ {\sl Let $a^*=\Delta=0,{\rm Re}(\mu)=0$ and 
$\alpha<1-\nu<\beta.$
\par
{\bf (a)} \ The transform {\boldmath${H}$} is defined on $\LL_{\nu,2}$ and it 
can be extended to $\LL_{\nu,r}$ as an element of $[\LL_{\nu,r},\LL_{1-\nu,r}]$ 
for $1<r<\infty$.
\par
{\bf (b)} \ If $1<r\eql2,$ the transform {\boldmath${H}$} is one-to-one on 
$\LL_{\nu,r}$ and there holds the equality
\begin{eqnarray}
&&\left(\M\mbox{\boldmath$H$}f\right)(s)=\HH(s)\left(\M f\right)(1-s)\quad
({\rm Re}(s)=1-\nu).
\end{eqnarray}
\par
{\bf (c)} \ If $\nu\notin{\cal E}_{\HHs},$ then {\boldmath${H}$} is one-to-one 
on $\LL_{\nu,r}$ and there holds
\begin{eqnarray}
&&\mbox{\boldmath$H$}(\LL_{\nu,r})=\LL_{1-\nu,r}.
\end{eqnarray}
\par
{\bf (d)} \ If $f\in \LL_{\nu,r}$ and $g\in \LL_{\nu,r'}$ with $1<r<\infty$ and 
$r'=r/(r-1),$ then the relation $(3.7)$ holds.
\par
{\bf (e)} \ If $f\in \LL_{\nu,r}$ with $1<r<\infty$ and $\lambda\in\C,h>0,$ 
then {\boldmath${H}$}$f$ is given by $(3.8)$ for 
${\rm Re}(\lambda)>(1-\nu)h-1,$ while {\boldmath${H}$}$f$ is given by $(3.9)$ 
for} ${\rm Re}(\lambda)<(1-\nu)h-1$.
\par
{\bf Proof.} \ Since $\alpha<1-\nu<\beta$ and $\Delta(1-\nu)+{\rm Re}(\mu)
\eql0,$ then according to Theorem 3 the transform {\boldmath${H}$} is defined 
on $f\in \LL_{\nu,2}$. We denote by $\HH_0(s)$ the function
\begin{eqnarray}
&&\HH_0(s)=\delta^{-s}\HH(s),
\end{eqnarray}
where $\delta$ is defined in (1.10). It follows from (3.3) that
\begin{eqnarray}
&&\left|\HH_0(\sigma+it)\right|\sim\prod^p_{i=1}\alpha_i^{1/2-{\rm Re}(a_i)}
\prod^q_{j=1}\beta_j^{{\rm Re}(b_j)-1/2}(2\pi)^{c^*}e^{-\pi{\rm Im}(\xi)
{\rm sign}(t)/2}\quad(|t|\to\infty)
\end{eqnarray}
is uniformly in $\sigma$ for $\sigma$ in any bounded interval in $\R$. 
Therefore $\HH_0(s)$ is analytic in the strip $\alpha<{\rm Re}(s)<\beta$, and 
if $\alpha<\sigma_1\eql\sigma_2<\beta$, then $|\HH_0(s)|$ is bounded in the 
strip $\sigma_1\eql{\rm Re}(s)\eql\sigma_2$. Since $a^*=\Delta=0$, then in 
accordance with (1.9) $a^*_1=-a^*_2=\Delta/2=0$. Then from (4.3) and (3.5) we 
have
\begin{eqnarray}
\HH'_0(\sigma+it)&\hspace{-2.5mm}=&\hspace{-2.5mm}\HH_0(\sigma+it)
\left[-\log(\delta)+\frac{\HH'(\sigma+it)}{\HH(\sigma+it)}\right]\\[2mm]
&\hspace{-2.5mm}=&\hspace{-2.5mm}\HH_0(\sigma+it)\left[-\log(\delta)+
\log(\delta)+\frac{{\rm Im}(\mu)}{it}+O\left(\frac1{t^2}\right)\right]\nonumber
\\[2mm]
&\hspace{-2.5mm}=&\hspace{-2.5mm}O\left(\frac1t\right)\quad(|t|\to\infty)
\nonumber
\end{eqnarray}
for $\alpha<\sigma<\beta$. Thus $\HH_0(s)$ belongs to the class ${\cal A}$ 
(see Definition 1) with $\alpha(\HH_0)=\alpha$ and $\beta(\HH_0)=\beta$. 
Therefore by virtue of Theorem 1, there is a transform $T\in[\LL_{\nu,r}]$ 
with $1<r<\infty$ and $\alpha<\nu<\beta$. When $1<r\eql2$, then $T$ is 
one-to-one on $\LL_{\nu,2}$ and the relation
\begin{eqnarray}
&&\left(\M Tf\right)(s)=\HH_0(s)\left(\M f\right)(s)\quad({\rm Re}(s)=\nu)
\end{eqnarray}
holds for $f\in \LL_{\nu,r}$. Let
\begin{eqnarray}
&&\mbox{\boldmath$H$}_0=W_\delta TR,
\end{eqnarray}
where $W_\delta$ and $R$ are given by (2.14) and (2.15). According to the 
properties (P2) and (P3) of the operators $W_\delta$ and $R$, we find 
$R\in[\LL_{\nu,r},\LL_{1-\nu,r}]$, $W_\delta\in[\LL_{1-\nu,r},\LL_{1-\nu,r}]$ 
and hence {\boldmath$H$}$_0\in[\LL_{\nu,r},\LL_{1-\nu,r}]$ for 
$\alpha<1-\nu<\beta$ and $1<r<\infty$, too. When $\alpha<1-\nu<\beta$, 
$1<r\eql2$ and $f\in \LL_{\nu,r}$, it follows from (4.7), (2.17), (4.6), (2.18) 
and (4.3) that
\begin{eqnarray}
\left(\M\mbox{\boldmath$H$}_0f\right)(s)&\hspace{-2.5mm}=&\hspace{-2.5mm}
\left(\M W_\delta TRf\right)(s)=\delta^s\left(\M TRf\right)(s)\\[2mm]
&\hspace{-2.5mm}=&\hspace{-2.5mm}\delta^s\HH_0(s)\left(\M Rf\right)(s)=
\delta^s\HH_0(s)\left(\M f\right)(1-s)=\HH(s)(\M f)(1-s)\nonumber
\end{eqnarray}
for ${\rm Re}(s)=1-\nu$. In particular, for $f\in \LL_{\nu,2}$ Theorem 3 (i), 
(3.6) and (4.8) imply the equality
\begin{eqnarray}
&&\left(\M\mbox{\boldmath$H$}_0f\right)(s)=\left(\M\mbox{\boldmath$H$}f\right)
(s)\quad({\rm Re}(s)=1-\nu).
\end{eqnarray}
Thus {\boldmath$H$}$_0f=${\boldmath$H$}$f$ for $f\in \LL_{\nu,2}$ and 
therefore if $\alpha<1-\nu<\beta$, {\boldmath$H$}$=${\boldmath$H$}$_0$ on 
$\LL_{\nu,2}$ by Theorem 3 (iv). Since $\LL_{\nu,2}\bigcap \LL_{\nu,r}$ is 
dense in $\LL_{\nu,r}$ [39, Theorem 2.2], {\boldmath$H$} can be extended to 
$\LL_{\nu,r}$ if we define it there by {\boldmath$H$}$_0$, and then 
{\boldmath$H$}$\in[\LL_{\nu,r},\LL_{1-\nu,r}]$. This completes the proof of 
assertion (a) of the theorem.
\par
The assertions (b) - (e) are proved similarly to those in [19, Theorem 4.1] 
where the case $\delta=1$ was considered, if we take into account that 
$W_\delta$ is a one-to-one transform on $\LL_{1-\nu,r}$ and 
$W_\delta\left(\LL_{1-\nu,r}\right)=\LL_{1-\nu,r}$ for $\alpha<1-\nu<\beta$ 
and $1<r<\infty$. The theorem is proved.
\\\par
{\bf Theorem 6.} \ {\sl Let $a^*=\Delta=0,{\rm Re}(\mu)<0$ and 
$\alpha<1-\nu<\beta,$ and let either $m>0$ or $n>0$.
\par
{\bf (a)} \ The transform {\boldmath$H$} defined on $\LL_{\nu,2}$ can be 
extended to $\LL_{\nu,r}$ for $1<r<\infty$ as an element of 
$[\LL_{\nu,r},\LL_{1-\nu,s}]$ for all $s\eqg r$ such that 
$1/s>1/r+{\rm Re}(\mu)$.
\par
{\bf (b)} \ If $1<r\eql2,$ then {\boldmath$H$} is a one-to-one transform on 
$\LL_{\nu,r}$ and there holds the equality $(4.1).$
\par
{\bf (c)} \ If $\nu\notin{\cal E}_{\HHs},$ then {\boldmath$H$} is a one-to-one 
transform on $\LL_{\nu,r}$ and there hold
\begin{eqnarray}
&&\mbox{\boldmath$H$}\left(\LL_{\nu,r}\right)=I^{-\mu}_{-;k,-\alpha/k}
\left(\LL_{1-\nu,r}\right)
\end{eqnarray}
for $k\eqg1$ and $m>0,$ and
\begin{eqnarray}
&&\mbox{\boldmath$H$}\left(\LL_{\nu,r}\right)=I^{-\mu}_{0+;k,\beta/k-1}
\left(\LL_{1-\nu,r}\right)
\end{eqnarray}
for $0<k\eql1$ and $n>0$. If $\nu\in{\cal E}_{\HHs},$ 
{\boldmath$H$}$\left(\LL_{\nu,r}\right)$ is a subset of 
$I^{-\mu}_{-;k,-\alpha/k}\left(\LL_{\nu,r}\right)$ or 
$I^{-\mu}_{0+;k,\beta/k-1}\left(\LL_{\nu,r}\right),$ when $m>0$ or $n>0,$ 
respectively.
\par
{\bf (d)} \ If $f\in \LL_{\nu,r}$ and $g\in \LL_{\nu,s}$ with 
$1<r<\infty,1<s<\infty$ and $1\eql1/r+1/s<1-{\rm Re}(\mu),$ then the relation 
$(3.7)$ holds.
\par
{\bf (e)} \ If $f\in \LL_{\nu,r}$ with $1<r<\infty$ and $\lambda\in\C,h>0,$ 
then {\boldmath$H$}$f$ is given by $(1.8)$ for ${\rm Re}(\lambda)>(1-\nu)h-1,$ 
while {\boldmath$H$}$f$ is given by $(3.9)$ for ${\rm Re}(\lambda)<(1-\nu)h-1$. 
If furthermore ${\rm Re}(\mu)<-1,$ then {\boldmath$H$}$f$ is given by} (1.1).
\par
{\bf Proof.} \ Since $\alpha<1-\nu<\beta$, $a^*=0$ and 
$\Delta(1-\nu)+{\rm Re}(\mu)={\rm Re}(\mu)<0$, then from Theorem 3 the 
transform {\boldmath$H$} is defined on $\LL_{\nu,2}$.
\par
If $m>0$ or $n>0$, then $\alpha$ or $\beta$ are finite in view of (3.1) and 
(3.2). We set
\begin{eqnarray}
\HH_1(s)&\hspace{-2.5mm}=&\hspace{-2.5mm}\frac{\Gamma([s-\alpha]/k-\mu)}
{\Gamma([s-\alpha]/k)}\HH(s)\\[2mm]
&\hspace{-2.5mm}=&\hspace{-2.5mm}\HH^{m+1,n}_{\thinspace p+1,q+1}
\left[\left.\begin{array}{l}(a_1,\alpha_1),\cdots,(a_p,\alpha_p),
(-\alpha/k,1/k)\\[2mm](-\mu-\alpha/k,1/k),(b_1,\beta_1),\cdots,(b_q,\beta_q)
\end{array}\right|s\right]\nonumber
\end{eqnarray}
for $m>0$ and $k\eqg1$, and
\begin{eqnarray}
\HH_2(s)&\hspace{-2.5mm}=&\hspace{-2.5mm}\frac{\Gamma([\beta-s]/k-\mu)}
{\Gamma([\beta-s]/k)}\HH(s)\\[2mm]
&\hspace{-2.5mm}=&\hspace{-2.5mm}\HH^{m,n+1}_{\thinspace p+1,q+1}
\left[\left.\begin{array}{l}(1+\mu-\beta/k,1/k),(a_1,\alpha_1),
\cdots,(a_p,\alpha_p)\\[2mm](b_1,\beta_1),\cdots,(b_q,\beta_q),
(1-\beta/k,1/k),\end{array}\right|s\right]\nonumber
\end{eqnarray}
for $n>0$ and $0<k\eql1$. We denote by $\alpha_1,$ $\beta_1,$ 
$\widetilde{a}^*_1,$ $\Delta_1,$ $\delta_1$ and $\mu_1$ for $\HH_1$, and by 
$\alpha_2,$ $\beta_2,$ $\widetilde{a}^*_2,$ $\Delta_2,$ $\delta_2$ and $\mu_2$ 
for $\HH_2$ instead of that for \HH. Then we find that
\begin{eqnarray*}
\begin{array}{l}
\alpha_1=\max[\alpha,\alpha+k{\rm Re}(\mu)]=\alpha,\beta_1=\beta,
\widetilde{a}^*_1=a^*=0,\Delta_1=\Delta=0,\delta_1=\delta,\mu_1=0,\\[2mm]
\alpha_2=\alpha,\beta_2=\min[\beta,\beta-k{\rm Re}(\mu)]=\beta,
\widetilde{a}^*_2=a^*=0,\Delta_2=\Delta=0,\delta_2=\delta,\mu_2=0,
\end{array}
\end{eqnarray*}
and the exceptional sets ${\cal E}_{\HHs_1}$ and ${\cal E}_{\HHs_2}$ of 
$\HH_1$ and $\HH_2$ coincide with that ${\cal E}_{\HHs}$ of \HH. Then 
according to Theorem 5, if $\alpha<1-\nu<\beta$ and $1<r<\infty,$ there are 
transforms $\widetilde{\mbox{\boldmath$H$}}_1\in[\LL_{\nu,r},\LL_{1-\nu,r}]$ 
for $m>0$, $k\eqg1$ and $\widetilde{\mbox{\boldmath$H$}}_2\in[\LL_{\nu,r},
\LL_{1-\nu,r}]$ for $n>0$, $0<k\eql1$; and if $f\in \LL_{\nu,r}$ with 
$1<r\eql2$, then by (4.1)
\begin{eqnarray}
&&\left(\M\widetilde{\mbox{\boldmath$H$}}_i f\right)(s)=
\HH_i(s)\left(\M f\right)(1-s)\quad({\rm Re}(s)=1-\nu;\ i=1,2).
\end{eqnarray}
We set
\begin{eqnarray}
&&{\mbox{\boldmath$H$}}_1=I^{-\mu}_{-;k,-\alpha/k}
\widetilde{\mbox{\boldmath$H$}}_1\quad\mbox{for}\quad m>0,k\eqg1,
\end{eqnarray}
and
\begin{eqnarray}
&&{\mbox{\boldmath$H$}}_2=I^{-\mu}_{0+;k,\beta/k-1}
\widetilde{\mbox{\boldmath$H$}}_2\quad\mbox{for}\quad n>0,0<k\eql1.
\end{eqnarray}
Using Theorem 2 (b) in the first case and Theorem 2 (a) in the second, and 
taking the same arguments as in the proof of Theorem 5, we arrive at the 
assertion (a) of the theorem. The assertions (b) - (e) are proved similarly to 
that in [19, Theorem 4.2] (while considering the case $\delta=1$) on the basis 
of Theorem 5 and Theorem 4 in the case ${\rm Re}(\mu)<-1$.
\\\\
\begin{flushleft}
{\bf 5. \ $\LL_{\nu,r}$-Theory of the Transform {\boldmath$H$} 
($a^*=0,\Delta\ne0$)}
\end{flushleft}
\setcounter{section}{5}
\setcounter{equation}{0}
\par
In this section we discuss that, if $a^*=0$ and $\Delta\ne0$, then the 
transform {\boldmath$H$} defined on $\LL_{\nu,2}$ can be extended to 
$\LL_{\nu,r}$ such as $\mbox{\boldmath$H$}\in[\LL_{\nu,r},\LL_{1-\nu,s}]$ for 
some range of values $s$. Then we characterize the range {\boldmath$H$} on 
$\LL_{\nu,r}$, except for its isolated values $\nu\in{\cal E}_{\HHs}$ in terms 
of the Hankel modified transform $\mbox{\boldmath$H$}_{k,\eta}$ and the 
elementary transform $M_\zeta$ given in (2.6) and (2.13). The result will be 
different in the cases $\Delta>0$ and $\Delta<0$. First we consider the case 
$\Delta>0$.
\\\par
{\bf Theorem 7.} \ {\sl Let $a^*=0,\Delta>0,-\infty<\alpha<1-\nu<\beta,
1<r<\infty$ and $\Delta(1-\nu)+{\rm Re}(\mu)\eql1/2-\gamma(r),$ where 
$\gamma(r)$ is defined in  $(2.10).$
\par
{\bf (a)} \ The transform {\boldmath$H$} defined on $\LL_{\nu,2}$ can be 
extended to $\LL_{\nu,r}$ to be an element of $[\LL_{\nu,r},\LL_{1-\nu,s}]$ 
for all $s$ with $r\eql s<\infty$ such that 
$s'\eqg[1/2-\Delta(1-\nu)-{\rm Re}(\mu)]^{-1}$ with $1/s+1/s'=1$.
\par
{\bf (b)} \ If $1<r\eql2,$ the transform {\boldmath$H$} is one-to-one on 
$\LL_{\nu,r}$ and there holds the equality
\begin{eqnarray}
&&(\M\mbox{\boldmath$H$}f)(s)=\HH(s)(\M f)(1-s)\quad({\rm Re}(s)=1-\nu).
\end{eqnarray}
\par
{\bf (c)} \ If $\nu\notin{\cal E}_{\HHs},$ then {\boldmath$H$} is one-to-one 
transform on $\LL_{\nu,r}$. If we set $\eta=-\Delta\alpha-\mu-1,$ then 
${\rm Re}(\eta)<-1$ and there holds
\begin{eqnarray}
&&\mbox{\boldmath$H$}(\LL_{\nu,r})=\left(M_{\mu/\Delta+1/2}H_{\Delta,\eta}
M_{\mu/\Delta+1/2}\right)(\LL_{\nu,r})=\left(M_{\mu/\Delta+1/2}H_{\Delta,\eta}
\right)(\LL_{\nu-{\rm Re}(\mu)/\Delta-1/2,r}).
\end{eqnarray}
When $\nu\in{\cal E}_{\HHs},$ then $\mbox{\boldmath$H$}(\LL_{\nu,r})$ is a 
subset of the right hand side of $(5.2).$
\par
{\bf (d)} \ If $f\in \LL_{\nu,r}$ and $g\in \LL_{\nu,s}$ with 
$1<r<\infty,1<s<\infty,1/r+1/s\eqg1$ and $\Delta(1-\nu)+{\rm Re}(\mu)
\eql1/2-\max[\gamma(r),\gamma(s)],$ then the relation $(3.7)$ holds.
\par
{\bf (e)} \ If $f\in \LL_{\nu,r}$ with $1<r<\infty,\lambda\in\C,h>0$ and 
$\Delta(1-\nu)+{\rm Re}(\mu)\eql1/2-\gamma(r),$ then {\boldmath$H$}$f$ is 
given by $(3.8)$ for ${\rm Re}(\lambda)>(1-\nu)h-1,$ while {\boldmath$H$}$f$ 
is given by $(3.9)$ for ${\rm Re}(\lambda)<(1-\nu)h-1$. If 
$\Delta(1-\nu)+{\rm Re}(\mu)<-1,$ {\boldmath$H$}$f$ is given by} (1.1).
\par
{\bf Proof.} \ Since $\gamma(r)\eqg1/2$, we have $\Delta(1-\nu)+{\rm Re}(\mu)
\eql0$ by the assumption, and hence from Theorem 3 the transform 
{\boldmath$H$} is defined on $\LL_{\nu,2}$. The condition $\Delta>0$ and the 
relation $\Delta(1-\nu)+{\rm Re}(\mu)\eql0$ imply 
$\nu\eqg1+{\rm Re}(\mu)/\Delta$ and $\alpha<-{\rm Re}(\mu)/\Delta$.
\par
Since $a^*=0$, then by (1.9)
\begin{eqnarray}
&&a^*_1=-a^*_2=\frac\Delta2>0.
\end{eqnarray}
We denote by $\HH_3(s)$ the function
\begin{eqnarray}
&&\HH_3(s)=\delta^{s-1}(a^*_1)^{(1-s)\Delta+\mu}\ \frac{\Gamma(-\mu+a^*_1
[s-1-\alpha])}{\Gamma(a^*_1[1-\alpha-s])}\ \HH(1-s).
\end{eqnarray}
As already known, the function $\HH(1-s)$ is analytic in the strip 
$1-\beta<{\rm Re}(s)<1-\alpha$ and the function 
$\Gamma(-\mu+a^*_1[s-1-\alpha])$ is analytic in the half-plane 
${\rm Re}(s)>\alpha+1+{\rm Re}(\mu)/a^*_1$$=\alpha+1+2{\rm Re}(\mu)/\Delta$. 
Since $\alpha<-{\rm Re}(\mu)/\Delta$, then 
$1-\alpha>\alpha+1+2{\rm Re}(\mu)/\Delta$. Therefore if we take 
$\alpha_1=\max[1-\beta,\alpha+1+2{\rm Re}(\mu)/\Delta]$ and 
$\beta_1=1-\alpha$, then $\alpha_1<\beta_1$ and $\HH_3(s)$ is analytic in the 
strip $\alpha_1<{\rm Re}(s)<\beta_1$.
\par
Setting $s=\sigma+it$ and a complex constant $k=c+id$, we have the behavior
\begin{eqnarray}
&&\Gamma(s+k)=\Gamma(c+\sigma+i[d+t])\sim\sqrt{2\pi}|t|^{c+\sigma-1/2}
e^{-\pi|t|/2-\pi d\thinspace{\rm sign}(t)/2}
\end{eqnarray}
as $|t|\to\infty$ (see [18, (2.12)]). Then by taking $a^*=0$, $\Delta=2a^*_1$ 
and (3.3) into account we have from (5.4) that
\begin{eqnarray}
\left|\HH_3(\sigma+it)\right|&\hspace{-2.5mm}\sim&\hspace{-2.5mm}(2\pi)^{c^*}
\prod^p_{i=1}\alpha_i^{1/2-{\rm Re}(a_i)}\prod^q_{j=1}
\beta_j^{{\rm Re}(b_j)-1/2}\left(a^*_1\right)^{(1-\sigma)\Delta+{\rm Re}(\mu)}
\\[2mm]
&&\times\ \left|a^*_1t\right|^{2a^*_1\sigma-2a^*_1-{\rm Re}(\mu)}
|t|^{(1-\sigma)\Delta+{\rm Re}(\mu)}e^{\pi[{\rm Im}(\mu)-{\rm Im}(\xi)]
{\rm sign}(t)/2}=\kappa\ne0\nonumber
\end{eqnarray}
as $|t|\to\infty$, uniformly in $\sigma$ for $\sigma$ in any bounded interval. 
Therefore, if $\alpha_1<\sigma_1\eql\sigma_2<\beta_1$, then $\HH_3(s)$ is 
bounded in $\sigma_1\eql{\rm Re}(s)\eql\sigma_2$.
\par
If $\alpha_1<\sigma<\beta_1$, then
\begin{eqnarray}
\HH'_3(\sigma+it)&\hspace{-2.5mm}=&\hspace{-2.5mm}\HH_3(\sigma+it)
\Biggl\{\log(\delta)-\Delta\log(a^*_1)+a^*_1\psi(a^*_1[s-\alpha-1]-\mu)\\[2mm]
&&+a^*_1\psi(a^*_1[1-\alpha-s])-\frac{\HH'(1-s)}{\HH(1-s)}\Biggr\},\nonumber
\end{eqnarray}
where $\psi$ is the psi-function $\psi(z)=\Gamma'(z)/\Gamma(z)$. Applying the 
estimate
\begin{eqnarray}
&&\psi(c+\sigma+it)=\log(it)+\frac{c+\sigma-1/2}{it}+O\left(\frac1{t^2}\right)
\quad(|t|\to\infty)
\end{eqnarray}
with $c\in\C$ (see [40, (3.9a)] and using (3.5) with $a^*=0$, (5.3), (5.6) and 
(5.8), we have from (5.7), as $|t|\to\infty$,
\begin{eqnarray}
&&\hspace{-10mm}\HH'_3(\sigma+it)\\[2mm]
&&\hspace{-5mm}=\HH_3(\sigma+it)\Biggl\{\log(\delta)-2a^*_1\log(a^*_1)+
a^*_1\left[\log(ia^*_1t)+\frac{a^*_1(\sigma-\alpha-1)-\mu-1/2}{ia^*_1t}\right]
\nonumber\\[2mm]
&&+a^*_1\left[\log(-ia^*_1t)-\frac{a^*_1(1-\alpha-\sigma)-1/2}{ia^*_1t}\right]
\nonumber\\[2mm]
&&-\left[\log(\delta)+a^*_1\log(-it)+a^*_1\log(it)-\frac{\mu+\Delta(1-\sigma)}
{it}\right]+O\left(\frac1{t^2}\right)\Biggr\}=O\left(\frac1{t^2}\right).
\nonumber
\end{eqnarray}
So $\HH_3\in{\cal A}$ with $\alpha(\HH_3)=\alpha_1$ and $\beta(\HH_3)=\beta_1$. 
Hence due to Theorem 1 there is a transform $T_3\in[\LL_{\nu,r}]$ corresponding 
to $\HH_3$ with $1<r<\infty$ and $\alpha_1<\nu<\beta_1$, and if $1<r\eql2$, 
then $T_3$ is a one-to-one on $\LL_{\nu,r}$ and
\begin{eqnarray}
&&(\M T_3f)(s)=\HH_3(s)(\M f)(s)\quad({\rm Re}(s)=\nu).
\end{eqnarray}
In particular, if $\nu$ and $r$ satisfy the hypothesis of this theorem, it is 
directly verified that $\alpha_1<\nu<\beta_1$ and hence the relation (5.10) is 
true.
\par
For $\eta=-\Delta\alpha-\mu-1$ let {\boldmath$H$}$_3$ be the operator
\begin{eqnarray}
&&\mbox{\boldmath$H$}_3=W_\delta M_{\mu/\Delta+1/2}H_{\Delta,\eta}
M_{\mu/\Delta+1/2}T_3
\end{eqnarray}
composed by the operator $W_\delta$ in (2.14), the operator $M_\zeta$ in 
(2.13), the modified Hankel transform (2.6) and the transform $T_3$ above, 
where $\alpha<-{\rm Re}(\mu)/\Delta$ so that ${\rm Re}(\eta)>-1$. For 
$1<r<\infty$ and $\alpha_1<\nu<\beta_1$ the properties (P1), (P2) and 
Theorem 2(c) yield {\boldmath$H$}$_3\in[\LL_{r,\nu},\LL_{1-\nu,s}]$ for all 
$s\eqg r$ such that $s'\eqg[1/2-\Delta(1-\nu)-{\rm Re}(\mu)]^{-1}$, and in 
particular, {\boldmath$H$}$_3\in[\LL_{r,\nu},\LL_{1-\nu,r}]$.
\par
If $f\in \LL_{\nu,r}$ with $1<r\eql2$, then applying (2.17), (2.16), (2.11), 
(5.10) and (5.4) and using the relation $\eta=-\Delta\alpha-\mu-1$, we have 
for ${\rm Re}(s)=1-\nu$
\begin{eqnarray}
&&\hspace{-3mm}\left(\M\mbox{\boldmath$H$}_3f\right)(s)\\[2mm]
&&=\left(\M W_\delta M_{\mu/\Delta+1/2}H_{\Delta,\eta}M_{\mu/\Delta+1/2}T_3f
\right)(s)\nonumber\\[2mm]
&&=\delta^s\left(\M M_{\mu/\Delta+1/2}H_{\Delta,\eta}M_{\mu/\Delta+1/2}T_3f
\right)(s)\nonumber\\[2mm]
&&=\delta^s\left(\M H_{\Delta,\eta}M_{\mu/\Delta+1/2}T_3f\right)
\left(s+\frac\mu\Delta+\frac12\right)\nonumber\\[2mm]
&&=\delta^s\left(\frac2\Delta\right)^{\Delta s+\mu}
\frac{\Gamma([\eta+\Delta s+\mu+1]/2)}{\Gamma([\eta-\Delta s-\mu+1]/2)}
\left(\M M_{\mu/\Delta+1/2}T_3f\right)\left(1-s-\frac\mu\Delta-\frac12\right)
\nonumber\\[2mm]
&&=\delta^s\left(\frac2\Delta\right)^{\Delta s+\mu}
\frac{\Gamma(\Delta[s-\alpha]/2)}{\Gamma(-\mu-\Delta[s+\alpha]/2)}\left(\M T_3
f\right)\left(1-s\right)\nonumber\\[2mm]
&&=\delta^s\left(a^*_1\right)^{-\mu-\Delta s}\frac{\Gamma(a^*_1[s-\alpha])}
{\Gamma(-\mu-a^*_1[s+\alpha])}\HH_3(1-s)\left(\M f\right)\left(1-s\right)
\nonumber\\[2mm]
&&=\HH(s)\left(\M f\right)\left(1-s\right).\nonumber
\end{eqnarray}
\par
In particular, if we take $r=2$ and $f\in \LL_{\nu,2}$, then 
$\left(\M\mbox{\boldmath$H$}_3f\right)(s)=\left(\M\mbox{\boldmath$H$}
f\right)(s)$ for ${\rm Re}(s)=1-\nu$, and hence 
{\boldmath$H$}$_3=${\boldmath$H$} on $\LL_{\nu,2}$. Thus, for all $\nu$ and 
$r$ satisfying the hypotheses of this theorem, {\boldmath$H$} can be extended 
from $\LL_{\nu,2}$ to $\LL_{\nu,r}$ if we define it by {\boldmath$H$}$_3$ 
given in (5.11) as an operator on $[\LL_{\nu,r},\LL_{1-\nu,s}]$. This 
completes the proof of the statement (a) of the theorem.
\par
The assertions (b) - (e) are proved similarly to those in [19, Theorem 5.1], 
when the case $\delta=1$ are considered, and the theorem is proved.
\\\par
{\bf Theorem 8.} \ {\sl Let $a^*=0,\Delta<0,\alpha<1-\nu<\beta<\infty,
1<r<\infty$ and $\Delta(1-\nu)+{\rm Re}(\mu)\eql1/2-\gamma(r).$
\par
{\bf (a)} \ The transform {\boldmath$H$} defined on $\LL_{\nu,2}$ can be 
extended to $\LL_{\nu,r}$ to be an element of $[\LL_{\nu,r},\LL_{1-\nu,r}]$ 
for all $s$ with $r\eql s<\infty$ such that 
$s'\eqg[1/2-\Delta(1-\nu)-{\rm Re}(\mu)]^{-1}$ with $1/s+1/s'=1$.
\par
{\bf (b)} \ If $1<r\eql2,$ then the transform {\boldmath$H$} is one-to-one on 
$\LL_{\nu,r}$ and there holds the equality $(5.1).$
\par
{\bf (c)} \ If $\nu\notin{\cal E}_{\HHs},$ then {\boldmath$H$} is a one-to-one 
transform on $\LL_{\nu,r}$. If we set $\eta=-\Delta\beta-\mu-1,$ then 
${\rm Re}(\eta)>-1$ and the relation $(5.2)$ holds$.$ When 
$\nu\in{\cal E}_{\HHs},$ then {\boldmath$H$}$(\LL_{\nu,r})$ is a subset of 
the right hand side of $(5.2).$
\par
{\bf (d)} \ If $f\in \LL_{\nu,r}$ and $g\in \LL_{\nu,s}$ with 
$1<r<\infty,1<s<\infty,1/r+1/s\eqg1$ and $\Delta(1-\nu)+
{\rm Re}(\mu)\eql1/2-\max[\gamma(r),\gamma(s)],$ then the relation $(3.7)$ 
holds.
\par
{\bf (e)} \ If $f\in \LL_{\nu,r}$ with $1<r<\infty,\lambda\in\C,h>0$ and 
$\Delta(1-\nu)+{\rm Re}(\mu)\eql1/2-\gamma(r),$ then {\boldmath$H$}$f$ is 
given by $(3.8)$ for ${\rm Re}(\lambda)>(1-\nu)h-1,$ while {\boldmath$H$}$f$ 
is given by $(3.9)$ for ${\rm Re}(\lambda)<(1-\nu)h-1$. If 
$\Delta(1-\nu)+{\rm Re}(\mu)<-1,$ then {\boldmath$H$}$f$ is given by} (1.1).
\\\par
The proof of this theorem is based on Theorem 7 and it is similarly to that in 
the case $\delta=1$ in [19, Theorem 5.2].
\\\par
{\bf Corollary.} \ {\sl Let $1<r<\infty,\alpha<\beta,a^*=0$ and one of the 
following conditions holds
\\[3mm]
$\begin{array}{lll}
{\bf (a)}&\Delta>0,&\alpha<\displaystyle\frac{1/2-{\rm Re}(\mu)-\gamma(r)}
\Delta;\\[2mm]
{\bf (b)}&\Delta<0,&\beta>\displaystyle\frac{1/2-{\rm Re}(\mu)-\gamma(r)}
\Delta;\\[2mm]
{\bf (c)}&\Delta=0,&{\rm Re}(\mu)\eql0.
\end{array}$
\\[3mm]
Then the transform {\boldmath$H$} can be defined on $\LL_{\nu,r}$ with} 
$\alpha<1-\nu<\beta$.
\\\\
\begin{flushleft}
{\bf 6. \ $\LL_{\nu,r}$-Theory of the Transform {\boldmath$H$} 
($a^*>0,a^*_1\eqg0$ and $a^*_2\eqg0$)}
\end{flushleft}
\setcounter{section}{6}
\setcounter{equation}{0}
\par
When $a^*>0$ and $\alpha<1-\nu<\beta$, then according to Theorem 4 the 
{\boldmath$H$}-transform is defined on $\LL_{\nu,2}$ and given by (1.1). Let 
us show that it can be extended to $\LL_{\nu,r}$ for any $1\eql r\eql\infty$. 
The next statement is proved similarly to that in case $\delta=1$ in [20, 
Theorem 2.3].
\\\par
{\bf Theorem 9.} \ {\sl Let $a^*>0,$ $\alpha<1-\nu<\beta$ and 
$1\eql r\eql s\eql\infty.$
\par
{\bf (a)} \ The {\boldmath$H$}-transform given in $(1.1)$ is defined on 
$\LL_{\nu,2}$ and can be extended to $\LL_{\nu,r}$ as an element of 
$[\LL_{\nu,r},\LL_{1-\nu,s}]$. If $1\eql r\eql2,$ then {\boldmath$H$} is a 
one-to-one transform from $\LL_{\nu,r}$ onto $\LL_{1-\nu,s}$.
\par
{\bf (b)} \ If $f\in \LL_{\nu,r}$ and $g\in \LL_{\nu,s'}$ with $1/s+1/s'=1,$ 
then the relation $(3.7)$ holds.}
\\\par
We give conditions for the transform {\boldmath$H$} to be one-to-one on 
$\LL_{\nu,r}$ and to characterize its range on $\LL_{\nu,r}$ except for its 
isolated values $\nu\in{\cal E}_{\HHs}$, in terms of the modified Laplace 
transform $L_{\kappa,\alpha}$ and the Erd\'elyi-Kober type fractional 
integration operators $I^\alpha_{0+;\sigma,\eta}$ and 
$I^\alpha_{-;\sigma,\eta}$ given in (2.7), (2.4) and (2.5). In this section we 
consider the case when $a^*>0,a^*_1\eqg0$ and $a^*_2\eqg0$. The results will 
be different along combinations of signs of $a^*_1$ and $a^*_2$. First we 
consider the case $a^*_1>0$ and $a^*_2>0$. 
\\\par
{\bf Theorem 10.} \ {\sl Let $a^*_1>0,$ $a^*_2>0,$ 
$-\infty<\alpha<1-\nu<\beta<\infty$ and $\omega=\mu+a^*_1\alpha-a^*_2\beta+1$ 
and let $1<r<\infty.$
\par
{\bf (a)} \ If $\nu\notin{\cal E}_{\HHs},$ or if $1\eql r\eql2,$ then 
{\boldmath$H$} is a one-to-one transform on $\LL_{\nu,r}$.
\par
{\bf (b)} \ If ${\rm Re}(\omega)\eqg0$ and $\nu\notin{\cal E}_{\HHs},$ then
\begin{eqnarray}
&&\mbox{\boldmath$H$}(\LL_{\nu,r})=\left(L_{a^*_1,\alpha}L_{a^*_2,
1-\beta-\omega/a^*_2}\right)(\LL_{1-\nu,r}).
\end{eqnarray}
When $\nu\in{\cal E}_{\HHs},$ $\mbox{\boldmath$H$}(\LL_{\nu,r})$ is a subset 
of the right hand side of $(6.1).$
\par
{\bf (c)} \ If ${\rm Re}(\omega)<0$ and $\nu\notin{\cal E}_{\HHs},$ then
\begin{eqnarray}
&&\mbox{\boldmath$H$}(\LL_{\nu,r})=\left(I^{-\omega}_{-;1/a^*_1,-a^*_1\alpha}
L_{a^*_1,\alpha}L_{a^*_2,1-\beta}\right)(\LL_{1-\nu,r}).
\end{eqnarray}
When $\nu\in{\cal E}_{\HHs},$ $\mbox{\boldmath$H$}(\LL_{\nu,r})$ is a subset 
of the right hand side of} (6.2).
\par
{\bf Proof.} \ We first consider the case ${\rm Re}(\omega)\eqg0$. We define 
$\HH_4(s)$ by
\begin{eqnarray}
&&\HH_4(s)=\frac{(a^*_1)^{a^*_1(s-\alpha)-1}(a^*_2)^{a^*_2(\beta-s)+\omega-1}}
{\Gamma(a^*_1[s-\alpha])\Gamma(a^*_2[\beta-s]+\omega)}\ \delta^{-s}\HH(s).
\end{eqnarray}
Since ${\rm Re}(\omega)\eqg0$, the function $\HH_4(s)$ is analytic in the 
strip $\alpha<{\rm Re}(s)<\beta$. According to (3.3) and (5.5) we have the 
estimate, as $|t|\to\infty$,
\begin{eqnarray}
\left|\HH_4(\sigma+it)\right|&\hspace{-2.5mm}\sim&\hspace{-2.5mm}
(a^*_1)^{a^*_1(\sigma-\alpha)-1}(a^*_2)^{a^*_2(\beta-\sigma)+
{\rm Re}(\omega)-1}\\[2mm]
&&\times\ \frac1{2\pi}(a^*_1|t|)^{a^*_1(\alpha-\sigma)+1/2}
(a^*_2|t|)^{a^*_2(\sigma-\beta)-{\rm Re}(\omega)+1/2}e^{[(a^*_1+a^*_2)
|t|-{\rm Im}(\omega){\rm sign}(t)]\pi/2}\nonumber\\[2mm]
&&\hspace{-10mm}\times\ \prod^p_{i=1}\alpha_i^{1/2-{\rm Re}(a_i)}\prod^q_{j=1}
\beta_j^{{\rm Re}(b_j)-1/2}(2\pi)^{c^*}|t|^{\Delta\sigma+{\rm Re}(\mu)}
e^{-[a^*|t|+{\rm Im}(\xi){\rm sign}(t)]\pi/2}\nonumber\\[2mm]
&\hspace{-2.5mm}\sim&\hspace{-2.5mm}\prod^p_{i=1}\alpha_i^{1/2-{\rm Re}(a_i)}
\prod^q_{j=1}\beta_j^{{\rm Re}(b_j)-1/2}(2\pi)^{c^*-1}(a^*_1a^*_2)^{-1/2}
e^{-[{\rm Im}(\xi)+{\rm Im}(\omega)]{\rm sign}(t)\pi/2}\nonumber
\end{eqnarray}
uniformly in $\sigma$ on any bounded interval in \R. Further, in accordance 
with (3.5) and (5.8)
\begin{eqnarray}
\HH_4'(\sigma+it)&\hspace{-2.5mm}=&\hspace{-2.5mm}\HH_4(\sigma+it)
\Biggl\{a^*_1\log(a^*_1)-a^*_2\log(a^*_2)-a^*_1\psi(a^*_1[s-\alpha])\\[2mm]
&&+a^*_2\psi(a^*_2[\beta-s]+\omega)-\log(\delta)+\frac{\HH'(\sigma+it)}
{\HH(\sigma+it)}\Biggr\}\nonumber\\[2mm]
&\hspace{-2.5mm}=&\hspace{-2.5mm}\HH_4(\sigma+it)
\Biggl\{a^*_1\log(a^*_1)-a^*_2\log(a^*_2)-a^*_1\left[\log(ia^*_1t)+
\frac{a^*_1(\sigma-\alpha)-1/2}{ia^*_1t}\right]\nonumber\\[2mm]
&&+a^*_2\left[\log(-ia^*_2t)-\frac{a^*_2(\beta-\sigma)+\omega-1/2}{ia^*_2t}
\right]-\log(\delta)\nonumber\\[2mm]
&&+\left[\log(\delta)+a^*_1\log(it)-a^*_2\log(-it)+\frac{\mu+\Delta\sigma}{it}
\right]+O\left(\frac1{t^2}\right)\Biggr\}\nonumber\\[2mm]
&\hspace{-2.5mm}=&\hspace{-2.5mm}O\left(\frac1{t^2}\right)\nonumber
\end{eqnarray}
as $|t|\to\infty$. So $\HH_4\in{\cal A}$ with $\alpha(\HH_4)=\alpha$ and 
$\beta(\HH_4)=\beta$ and Theorem 1 implies that there is a transform 
$T_4\in[\LL_{\nu,r}]$ for $1<r<\infty$ and $\alpha<\nu<\beta$ so that if 
$1<r\eql2$, then the relation
\begin{eqnarray}
(\M T_4f)(s)=\HH_4(s)(\M f)(s)\quad({\rm Re}(s)=\nu)
\end{eqnarray}
holds. Let
\begin{eqnarray}
\mbox{\boldmath$H$}_4=W_\delta L_{a^*_1,\alpha}L_{a^*_2,1-\beta-\omega/a^*_2}
T_4R,
\end{eqnarray}
where $W_\delta$ and $R$ are defined by (2.14) and (2.15). Then it follows 
from the properties (P2) and (P3) in Section 2 and Theorem 2 (d) that if 
$1<r\eql s<\infty$ and $\alpha<1-\nu<\beta$, then 
{\boldmath$H$}$_4\in[\LL_{\nu,r},\LL_{1-\nu,s}]$.
\par
For $f\in\LL_{\nu,2}$ applying (6.7), (2.17), (2.12), (6.6) and (2.18), we 
have
\begin{eqnarray}
\left(\M\mbox{\boldmath$H$}_4f\right)(s)&\hspace{-2.5mm}=&\hspace{-2.5mm}
\left(\M W_\delta L_{a^*_1,\alpha}L_{a^*_2,1-\beta-\omega/a^*_2}T_4Rf\right)(s)
\\[2mm]&\hspace{-2.5mm}
=&\hspace{-2.5mm}\delta^s\frac{\Gamma(a^*_1[s-\alpha])}
{(a^*_1)^{a^*_1(s-\alpha)-1}}\left(\M L_{a^*_2,1-\beta-\omega/a^*_2}T_4Rf
\right)(1-s)\nonumber\\[2mm]&\hspace{-2.5mm}
=&\hspace{-2.5mm}\delta^s\frac{\Gamma(a^*_1[s-\alpha])}
{(a^*_1)^{a^*_1(s-\alpha)-1}}\frac{\Gamma\{a^*_2[1+s-(1-\beta-\omega/a^*_2)]\}}
{(a^*_2)^{a^*_2[1-s-(1-\beta-\omega/a^*_2)]-1}}\left(\M T_4Rf\right)(s)
\nonumber\\[2mm]&\hspace{-2.5mm}
=&\hspace{-2.5mm}\delta^s\frac{\Gamma(a^*_1[s-\alpha])}
{(a^*_1)^{a^*_1(s-\alpha)-1}}\frac{\Gamma(a^*_2[\beta-s]+\omega)}
{(a^*_2)^{a^*_2(\beta-s)+\omega-1}}\ \HH_4(s)\left(\M Rf\right)(s)\nonumber
\\[2mm]
&\hspace{-2.5mm}=&\hspace{-2.5mm}\HH(s)\left(\M f\right)(1-s).\nonumber
\end{eqnarray}
Thus we obtain that if $f\in \LL_{\nu,2}$, then 
$\left(\M\mbox{\boldmath$H$}_4f\right)(s)=
\left(\M\mbox{\boldmath$H$}f\right)(s)$ with ${\rm Re}(s)=1-\nu$. Hence 
{\boldmath$H$}$_4=${\boldmath$H$} on $\LL_{\nu,2}$ and {\boldmath$H$} can be 
extended from $\LL_{\nu,2}$ to $\LL_{\nu,r}$ if we define it by (6.7). Further 
the proof of assertion (b) is carried out similarly to that in [20, 
Theorem 3.1] for the case $\delta=1$.
\par
Let ${\rm Re}(\omega)<0$. We denote by $\HH_5(s)$ the function
\begin{eqnarray}
\HH_5(s)=\frac{(a^*_1)^{a^*_1(s-\alpha)-1}(a^*_2)^{a^*_2(\beta-s)-1}
\Gamma(a^*_1[s-\alpha]-\omega)}{\Gamma^2(a^*_1[s-\alpha])
\Gamma(a^*_2[\beta-s])}\ \delta^{-s}\HH(s),
\end{eqnarray}
which is analytic in the strip $\alpha<{\rm Re}(s)<\beta$. Similar arguments 
to (6.4) and (6.5) show the estimates
\begin{eqnarray}
&&|\HH_5(\sigma+it)|\\[2mm]
&&\sim\prod^p_{i=1}\alpha_i^{1/2-{\rm Re}(a_i)}\prod^q_{j=1}
\beta_j^{{\rm Re}(b_j)-1/2}(2\pi)^{c^*-1}
(a^*_1)^{-{\rm Re}(\omega)-1/2}(a^*_2)^{-1/2}e^{[-{\rm Im}(\xi)+
{\rm Im}(\omega)]{\rm sign}(t)\pi/2}\nonumber
\end{eqnarray}
and
\begin{eqnarray}
\HH_5'(\sigma+it)&\hspace{-2.5mm}=&\hspace{-2.5mm}\HH_5(\sigma+it)
\Biggl\{a^*_1\log(a^*_1)-a^*_2\log(a^*_2)+a^*_1\psi(a^*_1[s-\alpha]-\omega)
\\[2mm]&&-2a^*_1\psi(a^*_1[s-\alpha])+a^*_2\psi(a^*_2[\beta-s])-\log(\delta)+
\frac{\HH'(\sigma+it)}{\HH(\sigma+it)}\Biggr\}\nonumber\\[2mm]
&\hspace{-2.5mm}=&\hspace{-2.5mm}O\left(\frac1{t^2}\right),\nonumber
\end{eqnarray}
as $|t|\to\infty$ hold uniformly in $\sigma$ on any bounded interval in \R. So 
$\HH_5\in{\cal A}$ with $\alpha(\HH_5)=\alpha$ and $\beta(\HH_5)=\beta$. By 
Theorem 1 there is a transform $T_5\in[\LL_{\nu,r}]$ for $1<r<\infty$ and 
$\alpha<\nu<\beta$ so that, if $1<r\eql2$,
\begin{eqnarray}
(\M T_5f)(s)=\HH_5(s)(\M f)(s)\quad({\rm Re}(s)=\nu).
\end{eqnarray}
Let
\begin{eqnarray}
\mbox{\boldmath$H$}_5=W_\delta I^{-\omega}_{-;1/a^*_1,a^*_1\alpha}
L_{a^*_1,\alpha}L_{a^*_2,1-\beta}T_5R.
\end{eqnarray}
Using again the properties (P2), (P3) in Section 2 and Theorem 2 (b), (d), we 
have that if $\alpha<1-\nu<\beta$ and $1<r\eql s<\infty$, then 
{\boldmath$H$}$_5\in[\LL_{\nu,r},\LL_{1-\nu,s}]$. For $f\in \LL_{\nu,2}$ and 
${\rm Re}(s)=1-\nu$ applying (6.13), (2.17), (2.9), (2.12), (6.12), (2.18) and 
(6.9) we obtain similarly to (6.8) that
\begin{eqnarray}
&&\hspace{-10mm}\left(\M\mbox{\boldmath$H$}_5f\right)(s)\\[2mm]
&\hspace{-2.5mm}=&\hspace{-2.5mm}\delta^s\left(\M I^{-\omega}_{-;1/a^*_1,
-a^*_1\alpha}L_{a^*_1,\alpha}L_{a^*_2,1-\beta}T_5Rf\right)(s)\nonumber\\[2mm]
&\hspace{-2.5mm}=&\hspace{-2.5mm}\delta^s\frac{\Gamma(a^*_1[s-\alpha])}
{\Gamma(a^*_1[s-\alpha]-\omega)}\left(\M L_{a^*_1,\alpha}L_{a^*_2,1-\beta}
T_5Rf\right)(s)\nonumber\\[2mm]
&\hspace{-2.5mm}=&\hspace{-2.5mm}\delta^s\frac{\Gamma(a^*_1[s-\alpha])}
{\Gamma(a^*_1[s-\alpha]-\omega)}\frac{\Gamma(a^*_1[s-\alpha[)}
{(a^*_1)^{a^*_1(s-\alpha)-1}}\left(\M L_{a^*_2,1-\beta}T_5Rf\right)(1-s)
\nonumber\\[2mm]
&\hspace{-2.5mm}=&\hspace{-2.5mm}\delta^s\frac{\Gamma(a^*_1[s-\alpha])}
{\Gamma(a^*_1[s-\alpha]-\omega)}\frac{\Gamma(a^*_1[s-\alpha])}
{(a^*_1)^{a^*_1(s-\alpha)-1}}\frac{\Gamma(a^*_2[1-s]-[1-\beta])}
{(a^*_2)^{a^*_2((1-s)-(1-\beta))-1}}\left(\M T_5Rf\right)(s)\nonumber\\[2mm]
&\hspace{-2.5mm}=&\hspace{-2.5mm}\delta^s\frac{\Gamma^2(a^*_1[s-\alpha])
\Gamma(a^*_2[\beta-s])}{\Gamma(a^*_1[s-\alpha]-\omega)
(a^*_1)^{a^*_1(s-\alpha)-1}(a^*_2)^{a^*_2(\beta-s)-1}}\ \HH_5(s)
\left(\M Rf\right)(s)\nonumber\\[2mm]
&\hspace{-2.5mm}=&\hspace{-2.5mm}\HH(s)(\M f)(1-s).\nonumber
\end{eqnarray}
Applying this equality and using similar arguments to those in the case 
${\rm Re}(\omega)\eqg0$, we complete the proof for ${\rm Re}(\omega)<0$.
\\\par
Now we proceed to the case $a^*_1>0$ and $a^*_2=0$.
\\\par
{\bf Theorem 11.} \ {\sl Let $a^*_1>0,$ $a^*_2=0,$ 
$-\infty<\alpha<1-\nu<\beta,$ $\omega=\mu+a^*_1\alpha+1/2$ and $1<r<\infty.$
\par
{\bf (a)} \ If $\nu\notin{\cal E}_{\HHs},$ or $1<r\eql2,$ then {\boldmath$H$} 
is a one-to-one transform on $\LL_{\nu,r}$.
\par
{\bf (b)} \ If ${\rm Re}(\omega)\eqg0$ and $\nu\notin{\cal E}_{\HHs},$ then
\begin{eqnarray}
\mbox{\boldmath$H$}(\LL_{\nu,r})=L_{a^*_1,\alpha-\omega/a^*_1}(\LL_{\nu,r}).
\end{eqnarray}
When $\nu\in{\cal E}_{\HHs},$ $\mbox{\boldmath$H$}(\LL_{\nu,r})$ is a subset 
of the right hand side of $(6.15).$
\par
{\bf (c)} \ If ${\rm Re}(\omega)<0$ and $\nu\notin{\cal E}_{\HHs},$ then
\begin{eqnarray}
\mbox{\boldmath$H$}(\LL_{\nu,r})=\left(I^{-\omega}_{-;1/a^*_1,-a^*_1\alpha}
L_{a^*_1,\alpha}\right)(\LL_{\nu,r}).
\end{eqnarray}
When $\nu\in{\cal E}_{\HHs},$ $\mbox{\boldmath$H$}(\LL_{\nu,r})$ is a subset 
of the right hand side of} (6.16).
\par
{\bf Proof.} \ We first consider the case ${\rm Re}(\omega)\eqg0$. We define 
$\HH_6(s)$ by
\begin{eqnarray}
&&\HH_6(s)=\frac{(a^*_1)^{a^*_1(s-\alpha)+\omega-1}}{\Gamma(a^*_1[s-\alpha]+
\omega)}\ \delta^{-s}\HH(s).
\end{eqnarray}
Since ${\rm Re}(\omega)\eqg0$, $\HH_6(s)$ is analytic in the strip 
$\alpha<{\rm Re}(s)<\beta$. Arguments similar to those in (6.4) and (6.5) lead 
to the estimates
\begin{eqnarray}
&&|\HH_6(\sigma+it)|\sim\prod^p_{i=1}\alpha_i^{1/2-{\rm Re}(a_i)}
\prod^q_{j=1}\beta_j^{{\rm Re}(b_j)-1/2}(2\pi)^{c^*-1/2}(a^*_1)^{-1/2}
e^{[-{\rm Im}(\xi)+{\rm Im}(\omega)]{\rm sign}(t)\pi/2},
\end{eqnarray}
and
\begin{eqnarray}
&&\hspace{-15mm}\HH_6'(\sigma+it)\\[2mm]
&\hspace{-2.5mm}=&\hspace{-2.5mm}\HH_6(\sigma+it)\left\{a^*_1\log(a^*_1)-a^*_1
\psi(a^*_1[s-\alpha]+\omega)-\log(\delta)+\frac{\HH'(\sigma+it)}
{\HH(\sigma+it)}\right\}\nonumber\\[2mm]
&\hspace{-2.5mm}=&\hspace{-2.5mm}O\left(\frac1{t^2}\right)\quad(|t|\to\infty)
\nonumber
\end{eqnarray}
uniformly in $\sigma$ in any bounded interval in \R. Thus $\HH_6\in{\cal A}$ 
with $\alpha(\HH_6)=\alpha$ and $\beta(\HH_6)=\beta$ and Theorem 1 implies 
that there is a transform $T_6\in[\LL_{\nu,r}]$ for $1<r<\infty$ and 
$\alpha<\nu<\beta$ and, if $1<r\eql2$,
\begin{eqnarray}
(\M T_6f)(s)=\HH_6(s)(\M f)(s)\quad({\rm Re}(s)=\nu).
\end{eqnarray}
We set
\begin{eqnarray}
\mbox{\boldmath$H$}_6=W_\delta L_{a^*_1,\alpha-\omega/a^*_1}RT_6R.
\end{eqnarray}
Then it follows from the properties (P2) and (P3) in Section 2 and Theorem 2 
(d) that, if $\alpha<\nu<\beta$ and $1<r\eql s<\infty$, then 
{\boldmath$H$}$_6\in[\LL_{\nu,r},\LL_{1-\nu,s}]$.
\par
For $f\in \LL_{\nu,2}$ applying (6.21), (2.17), (2.12), (6.20), (2.18) and 
(6.17), we obtain for ${\rm Re}(s)=1-\nu$
\begin{eqnarray}
(\M\mbox{\boldmath$H$}_6f)(s)&\hspace{-2.5mm}=&\hspace{-2.5mm}
(\M W_\delta L_{a^*_1,\alpha-\omega/a^*_1}RT_6Rf)(s)=
\delta^s(\M L_{a^*_1,\alpha-\omega/a^*_1}RT_6Rf)(s)\\[2mm]
&\hspace{-2.5mm}=&\hspace{-2.5mm}\delta^s\frac{\Gamma(a^*_1[s-\alpha+
\omega/a^*_1])}{(a^*_1)^{a^*_1(s-\alpha+\omega/a^*_1)-1}}\ (\M RT_6Rf)(1-s)
\nonumber\\[2mm]
&\hspace{-2.5mm}=&\hspace{-2.5mm}\delta^s\frac{\Gamma(a^*_1[s-\alpha]+
\omega)}{(a^*_1)^{a^*_1(s-\alpha)+\omega-1}}\ (\M T_6Rf)(s)\nonumber\\[2mm]
&\hspace{-2.5mm}=&\hspace{-2.5mm}\delta^s\frac{\Gamma(a^*_1[s-\alpha]+
\omega)}{(a^*_1)^{a^*_1(s-\alpha)+\omega-1}}\ \HH_6(s)(\M Rf)(s)\nonumber\\[2mm]
&\hspace{-2.5mm}=&\hspace{-2.5mm}\delta^s\frac{\Gamma(a^*_1[s-\alpha]+
\omega)}{(a^*_1)^{a^*_1(s-\alpha)+\omega-1}}\ \HH_6(s)(\M f)(1-s)=
\HH(s)(\M f)(1-s).\nonumber
\end{eqnarray}
Applying this relation and using the arguments similarly to those in the case 
${\rm Re}(\omega)\eqg0$ of Theorem 10, we complete the proof of theorem for 
${\rm Re}(\omega)\eqg0$.
\par
For the case ${\rm Re}(\omega)<0$. We define $\HH_7(s)$ by
\begin{eqnarray}
\HH_7(s)=\frac{(a^*_1)^{a^*_1(s-\alpha)-1}\Gamma(a^*_1[s-\alpha]-\omega)}
{\Gamma^2(a^*_1[s-\alpha])}\ \delta^{-s}\HH(s),
\end{eqnarray}
which is analytic in the strip $\alpha<{\rm Re}(s)<\beta$, and in accordance 
with (3.3), (3.5), (5.5) and (5.8), we find
\begin{eqnarray}
|\HH_7(\sigma+it)|&\hspace{-2.5mm}\sim&\hspace{-2.5mm}\prod^p_{i=1}
\alpha_i^{1/2-{\rm Re}(a_i)}\prod^q_{j=1}\beta_j^{{\rm Re}(b_j)-1/2}
(2\pi)^{c^*-1/2}\\[2mm]
&&\times\ (a^*_1)^{-{\rm Re}(\omega)+1/2}e^{-[{\rm Im}(\xi)-{\rm Im}(\omega)]
{\rm sign}(t)\pi/2},\nonumber
\end{eqnarray}
\begin{eqnarray}
\HH_7'(\sigma+it)&\hspace{-2.5mm}=&\hspace{-2.5mm}\HH_7(\sigma+it)
\Biggl[a^*_1\log(a^*_1)+a^*_1\psi(a^*_1[s-\alpha]-\omega)\\[2mm]
&&-2a^*_1\psi(a^*_1[s-\alpha])-\log(\delta)+\frac{\HH'(\sigma+it)}
{\HH(\sigma+it)}\Biggr]=O\left(\frac1{t^2}\right)\nonumber
\end{eqnarray}
as $|t|\to\infty$ uniformly in $\sigma$ on any bounded interval in \R. Then 
$\HH_7\in{\cal A}$ with $\alpha(\HH_7)=\alpha$ and $\beta(\HH_7)=\beta$. By 
Theorem 1, there is $T_7\in[\LL_{\nu,r}]$ for $1<r<\infty$ and 
$\alpha<\nu<\beta$ so that, if $1<r\eql2$,
\begin{eqnarray}
(\M T_7f)(s)=\HH_7(s)(\M f)(s)\quad({\rm Re}(s)=\nu).
\end{eqnarray}
Setting
\begin{eqnarray}
\mbox{\boldmath$H$}_7=W_\delta I^{-\omega}_{-;1/a^*_1,-a^*_1\alpha}
L_{a^*_1,\alpha}RT_7R
\end{eqnarray}
according to the properties (P2) and (P3) in Section 2 and Theorem 2 (b), (d), 
we have that if $\alpha<1-\nu<\beta$ and $1<r\eql s<\infty$, then 
{\boldmath$H$}$_7\in[\LL_{\nu,r},\LL_{1-\nu,s}]$. If $f\in \LL_{\nu,2}$ and 
${\rm Re}(s)=1-\nu$, applying (6.27), (2.17), (2.9), (2.12), (2.18) and (6.26) 
we obtain similarly to (6.8) that
\begin{eqnarray}
(\M\mbox{\boldmath$H$}_7f)(s)&\hspace{-2.5mm}=&\hspace{-2.5mm}
\left(\M W_\delta I^{-\omega}_{-;1/a^*_1,-a^*_1\alpha}L_{a^*_1,\alpha}
RT_7Rf\right)(s)\\[2mm]
&\hspace{-2.5mm}=&\hspace{-2.5mm}\delta^s\left(\M I^{-\omega}_{-;1/a^*_1,
-a^*_1\alpha}L_{a^*_1,\alpha}RT_7Rf\right)(s)\nonumber\\[2mm]
&\hspace{-2.5mm}=&\hspace{-2.5mm}\delta^s\frac{\Gamma(a^*_1[s-\alpha])}
{\Gamma(a^*_1[s-\alpha]-\omega)}\left(\M L_{a^*_1,\alpha}RT_7Rf\right)(s)
\nonumber\\[2mm]
&\hspace{-2.5mm}=&\hspace{-2.5mm}\delta^s\frac{\Gamma^2(a^*_1[s-\alpha])}
{\Gamma(a^*_1[s-\alpha]-\omega)(a^*_1)^{a^*_1(s-\alpha)-1}}\left(\M RT_7Rf
\right)(1-s)\nonumber\\[2mm]
&\hspace{-2.5mm}=&\hspace{-2.5mm}\delta^s\frac{\Gamma^2(a^*_1[s-\alpha])}
{\Gamma(a^*_1[s-\alpha]-\omega)(a^*_1)^{a^*_1(s-\alpha)-1}}\left(\M T_7Rf
\right)(s)\nonumber\\[2mm]
&\hspace{-2.5mm}=&\hspace{-2.5mm}\delta^s\frac{\Gamma^2(a^*_1[s-\alpha])}
{\Gamma(a^*_1[s-\alpha]-\omega)(a^*_1)^{a^*_1(s-\alpha)-1}}\HH_7(s)
\left(\M Rf\right)(s)\nonumber\\[2mm]
&\hspace{-2.5mm}=&\hspace{-2.5mm}\HH(s)(\M f)(1-s).\nonumber
\end{eqnarray}
Using this relation and the arguments similar to those in the case 
${\rm Re}(\omega)<0$ of Theorem 10, we complete the proof for 
${\rm Re}(\omega)<0$.
\par
In the case $a^*_1=0$ and $a^*_2>0$ the following statement is proved on the 
basis of Theorem 11 similarly to that in the case $\delta=1$ in [20, 
Theorem 3.3].
\\\par
{\bf Theorem 12.} \ {\sl Let $a^*_1=0,$ $a^*_2>0,$ $\alpha<1-\nu<\beta<\infty$ 
and $\omega=\mu-a^*_2\beta+1/2$ and let $1<r<\infty.$
\par
{\bf (a)} \ If $\nu\notin{\cal E}_{\HHs},$ or $1<r\eql2,$ then {\boldmath$H$} 
is a one-to-one transform on $\LL_{\nu,r}$.
\par
{\bf (b)} \ If ${\rm Re}(\omega)\eqg0$ and $\nu\notin{\cal E}_{\HHs},$ then
\begin{eqnarray}
\mbox{\boldmath$H$}(\LL_{\nu,r})=L_{-a^*_2,\beta+\omega/a^*_2}(\LL_{\nu,r}).
\end{eqnarray}
When $\nu\in{\cal E}_{\HHs},$ $\mbox{\boldmath$H$}(\LL_{\nu,r})$ is a subset 
of the right hand side of $(6.29).$
\par
{\bf (c)} \ If ${\rm Re}(\omega)<0$ and $\nu\notin{\cal E}_{\HHs},$ then
\begin{eqnarray}
\mbox{\boldmath$H$}(\LL_{\nu,r})=\left(I^{-\omega}_{0+;1/a^*_2,a^*_2\beta-1}
L_{-a^*_2,\beta}\right)(\LL_{\nu,r}).
\end{eqnarray}
When $\nu\in{\cal E}_{\HHs},$ $\mbox{\boldmath$H$}(\LL_{\nu,r})$ is a subset 
of the right hand side of} (6.30).
\\\\
\begin{flushleft}
{\bf 7. \ $\LL_{\nu,r}$-Theory of {\boldmath$H$}-transform ($a^*>0$ and 
$a^*_1<0$ or $a^*_2<0$)}
\end{flushleft}
\setcounter{section}{7}
\setcounter{equation}{0}
\par
In this section we give conditions to the transform {\boldmath$H$} to be 
one-to-one on $\LL_{\nu,r}$ and characterize its range on $\LL_{\nu,r}$ except 
for its isolated values $\nu\in{\cal E}_{\HHs}$ in terms of the modified Hankel 
transform $H_{\kappa,\eta}$ and modified Laplace transform $L_{\kappa,\alpha}$ 
given in (2.6) and (2.7). The results will be different in the cases 
$a^*_1>0,a^*_2<0$ and $a^*_1<0,a^*_2>0$. We first consider the former case.
\\\par
{\bf Theorem 13.} \ {\sl Let $a^*>0,$ $a^*_1>0,$ $a^*_2<0,$ 
$\alpha<1-\nu<\beta$ and $1<r<\infty.$
\par
{\bf (a)} \ If $\nu\notin{\cal E}_{\HHs},$ or if $1<r\eql2,$ then 
{\boldmath$H$} is a one-to-one transform on $\LL_{\nu,r}$.
\par
{\bf (b)} \ Let
\begin{eqnarray}
\omega=a^*\eta-\mu-\frac12,
\end{eqnarray}
where $\mu$ is given by $(1.7),$ and let $\eta$ be chosen as
\begin{eqnarray}
a^*{\rm Re}(\eta)&\hspace{-2.5mm}\eqg&\hspace{-2.5mm}\gamma(r)+
2a^*_2(\nu-1)+{\rm Re}(\mu),\\[2mm]
{\rm Re}(\eta)&\hspace{-2.5mm}>&\hspace{-2.5mm}\nu-1,
\end{eqnarray}
with $\gamma(r)$ being given by $(2.10),$ and let $\zeta$ be chosen such that
\begin{eqnarray}
{\rm Re}(\zeta)<1-\nu.
\end{eqnarray}
If $\nu\notin{\cal E}_{\HHs},$ then
\begin{eqnarray}
&&\mbox{\boldmath$H$}(\LL_{\nu,r})=\left(M_{1/2+\omega/(2a^*_2)}
H_{-2a^*_2,2a^*_2\zeta+\omega-1}L_{-a^*_2,1/2+\eta-\omega/(2a^*_2)}\right)
\left(\LL_{3/2+{\rm Re}(\omega)/(2a^*_2)-\nu,r}\right),
\end{eqnarray}
where $M_\zeta$ is given by $(2.13).$ If $\nu\in{\cal E}_{\HHs},$ then 
$\mbox{\boldmath$H$}(\LL_{\nu,r})$ is a subset of the right hand side of} 
(7.5).
\par
{\bf Proof.} \ We consider the function
\begin{eqnarray}
&&\HH_8(s)=\frac{(a^*)^{a^*(s+\eta)-1}|a^*_2|^{-2a^*_2s-\omega}
\Gamma(a^*_2[s+\zeta]+\omega)}{\Gamma(a^*[s+\eta])\Gamma(a^*_2[\zeta-s])}
\ \delta^{-s}\HH(s).
\end{eqnarray}
For ${\rm Re}(s)=1-\nu$ according to (7.1), (7.2), (7.4) and the relations 
$\gamma(r)\eqg1/2$ and $-a^*_2>0$, we have
\begin{eqnarray*}
{\rm Re}[a^*_2(s+\zeta)+\omega]&\hspace{-2.5mm}=&\hspace{-2.5mm}
a^*_2[1-\nu+{\rm Re}(\zeta)]+a^*{\rm Re}(\eta)-{\rm Re}(\mu)-\frac12\\[2mm]
&\hspace{-2.5mm}\eqg&\hspace{-2.5mm}a^*_2[1-\nu+{\rm Re}(\zeta)]+
[\gamma(r)+2a^*_2(\nu-1)+{\rm Re}(\mu)]-{\rm Re}(\mu)-\frac12\\[2mm]
&\hspace{-2.5mm}=&\hspace{-2.5mm}a^*_2[\nu-1+{\rm Re}(\zeta)]+
\gamma(r)-\frac12\eqg a^*_2[\nu-1+{\rm Re}(\zeta)]>0
\end{eqnarray*}
and hence the function $\HH_8(s)$ is analytic in the strip 
$\alpha<{\rm Re}(s)<\beta$. Applying (3.3), (3.5), (5.5) and (5.8) we obtain 
the estimates
\begin{eqnarray}
&&|\HH_8(\sigma+it)|\\[2mm]
&&\sim\prod^p_{i=1}\alpha_i^{1/2-{\rm Re}(a_i)}\prod^q_{j=1}
\beta_j^{{\rm Re}(b_j)-1/2}(2\pi)^{c^*-1/2}(a^*)^{-1/2}
e^{-[2a^*_2{\rm Im}(\zeta)+{\rm Im}(\xi)+{\rm Im}(\omega)-a^*{\rm Im}(\eta)]
{\rm sign}(t)\pi/2}\nonumber
\end{eqnarray}
and
\begin{eqnarray}
\HH_8'(\sigma+it)&\hspace{-2.5mm}=&\hspace{-2.5mm}\HH_8(\sigma+it)
\Biggl\{a^*\log(a^*)-2a^*_2\log(|a^*_2|)+a^*_2\psi(a^*_2[s+\zeta]+\omega)\\[2mm]
&&-a^*\psi(a^*[s+\eta])+a^*_2\psi(a^*_2[\zeta-s])-\log(\delta)+
\frac{\HH'(\sigma+it)}{\HH(\sigma+it)}\Biggr\}\nonumber\\[2mm]
&\hspace{-2.5mm}=&\hspace{-2.5mm}O\left(\frac1{t^2}\right)\nonumber
\end{eqnarray}
as $|t|\to\infty$, uniformly in $\sigma$ on any bounded interval in \R.
\par
Thus we have $\HH_8\in{\cal A}$ with $\alpha(\HH_8)=\alpha$ and 
$\beta(\HH_8)=\beta$ and by Theorem 1, there is a transform 
$T_8\in[\LL_{\nu,r}]$ for $1<r<\infty$ and $\alpha<\nu<\beta$, so that if 
$1<r\eql2$,
\begin{eqnarray}
(\M T_8f)(s)=\HH_8(s)(\M f)(s)\quad({\rm Re}(s)=\nu).
\end{eqnarray}
Let
\begin{eqnarray}
&&\mbox{\boldmath$H$}_8=W_\delta M_{1/2+\omega/(2a^*_2)}
H_{-2a^*_2,2a^*_2\zeta+\omega-1}L_{-a^*,1/2+\eta-\omega/(2a^*_2)}
M_{-1/2-\omega/(2a^*_2)}T_8R.
\end{eqnarray}
It is directly verified that if the conditions of theorem satisfy, then in 
accordance with properties (P1) - (P3) in Section 2 and Theorem 2 (c), (d), 
{\boldmath$H$}$_8\in[\LL_{\nu,r},\LL_{1-\nu,r}]$.
\par
If $f\in \LL_{\nu,2}$, then applying (7.10), (2.17), (2.16), (2.11), (2.12), 
(2.16), (7.9), (2.18) and (7.6), we have for ${\rm Re}(s)=1-\nu$
\begin{eqnarray}
&&\hspace{-5mm}\left(\M\mbox{\boldmath$H$}_8f\right)(s)\\[2mm]
&&=\left(\M W_\delta M_{1/2+\omega/(2a^*_2)}H_{-2a^*_2,2a^*_2\zeta+\omega-1}
L_{-a^*,1/2+\eta-\omega/(2a^*_2)}M_{-1/2-\omega/(2a^*_2)}T_8Rf\right)(s)
\nonumber\\[2mm]
&&=\delta^s\left(\M M_{1/2+\omega/(2a^*_2)}H_{-2a^*_2,2a^*_2\zeta+\omega-1}
L_{-a^*,1/2+\eta-\omega/(2a^*_2)}M_{-1/2-\omega/(2a^*_2)}T_8Rf\right)(s)
\nonumber\\[2mm]
&&=\delta^s\left(\M H_{-2a^*_2,2a^*_2\zeta+\omega-1}L_{-a^*,1/2+
\eta-\omega/(2a^*_2)}M_{-1/2-\omega/(2a^*_2)}T_8Rf\right)\left(s+\frac12+
\frac\omega{2a^*_2}\right)\nonumber\\[2mm]
&&=\delta^s|a^*_2|^{2a^*_2s+\omega}\frac{\Gamma(a^*_2[\zeta-s])}
{\Gamma(a^*_2[s+\zeta]+\omega)}\nonumber\\[2mm]
&&\hspace{10mm}\times\ \left(\M L_{-a^*,1/2+\eta-\omega/(2a^*_2)}
M_{-1/2-\omega/(2a^*_2)}T_8Rf\right)\left(\frac12-s-\frac\omega{2a^*_2}\right)
\nonumber\\[2mm]
&&=\delta^s|a^*_2|^{2a^*_2s+\omega}(a^*)^{1-a^*(s+\eta)}
\frac{\Gamma(a^*_2[\zeta-s])\Gamma(a^*[s+\eta])}{\Gamma(a^*_2[s+\zeta]+\omega)}
\nonumber\\[2mm]
&&\hspace{10mm}\times\ \left(\M M_{-1/2-\omega/(2a^*_2)}T_8Rf\right)
\left(\frac12+s+\frac\omega{2a^*_2}\right)\nonumber\\[2mm]
&&=\delta^s|a^*_2|^{2a^*_2s+\omega}(a^*)^{1-a^*(s+\eta)}
\frac{\Gamma(a^*_2[\zeta-s])\Gamma(a^*[s+\eta])}
{\Gamma(a^*_2[s+\zeta]+\omega)}\left(\M T_8Rf\right)(s)\nonumber\\[2mm]
&&=\delta^s|a^*_2|^{2a^*_2s+\omega}(a^*)^{1-a^*(s+\eta)}
\frac{\Gamma(a^*_2[\zeta-s])\Gamma(a^*[s+\eta])}
{\Gamma(a^*_2[s+\zeta]+\omega)}\ \HH_8(s)\left(\M Rf\right)(s)\nonumber\\[2mm]
&&=\HH(s)\left(\M f\right)(1-s).\nonumber
\end{eqnarray}
Using this relation and the arguments similar to those in the case $\delta=1$ 
in [20, Theorem 5.1], we complete the proof of theorem.
\\\par
{\bf Corollary 1.} \ {\sl Let $a^*>0,$ $a^*_1>0,$ $a^*_2<0$ and let 
$-\infty\eql\alpha<1-\nu<\beta\eql\infty,$ $1<r<\infty.$
\par
{\bf (a)} \ If $\nu\notin{\cal E}_{\HHs},$ or if $1<r\eql2,$ then 
{\boldmath$H$} is a one-to-one transform on $\LL_{\nu,r}$.
\par
{\bf (b)} \ Let $\omega$ be given by $(7.1)$ and let $\eta$ and $\zeta$ be 
chosen such that either of the following conditions holds
\begin{enumerate}
\item[{{\bf (i)}\ }]$a^*{\rm Re}(\eta)\eqg\gamma(r)-2a^*_2\beta+{\rm Re}(\mu),$ 
${\rm Re}(\eta)\eqg-\alpha,$ ${\rm Re}(\zeta)\eql\alpha,$ if 
$-\infty<\alpha<\beta<\infty;$
\item[{{\bf (ii)}\ }]$a^*{\rm Re}(\eta)\eqg\gamma(r)-2a^*_2\beta
+{\rm Re}(\mu),$ ${\rm Re}(\eta)>\nu-1,$ ${\rm Re}(\zeta)<1-\nu,$ if 
$-\infty=\alpha<\beta<\infty;$
\item[{{\bf (iii)}\ }]$a^*{\rm Re}(\eta)\eqg\gamma(r)+2a^*_2(\nu-1)+
{\rm Re}(\mu),$ ${\rm Re}(\eta)\eqg-\alpha,$ ${\rm Re}(\zeta)\eql\alpha,$ if 
$-\infty<\alpha<\beta=\infty.$
\end{enumerate}
Then if $\nu\notin{\cal E}_{\HHs},$ {\boldmath$H$}$(\LL_{\nu,r})$ can be 
represented by the relation $(7.5),$ and if $\nu\in{\cal E}_{\HHs},$ 
{\boldmath$H$}$(\LL_{\nu,r})$ is a subset of the right hand side of} (7.5).
\\\par
{\bf Corollary 2.} \ {\sl Let $a^*>0,$ $a^*_1>0,$ $a^*_2<0$ and let 
$-\infty<\alpha<1-\nu<\beta<\infty,$ $1<r<\infty.$
\par
{\bf (a)} \ If $\nu\notin{\cal E}_{\HHs},$ or if $1<r\eql2,$ then 
{\boldmath$H$} is a one-to-one transform on $\LL_{\nu,r}$.
\par
{\bf (b)} \ Let $a^*\alpha-2a^*_2\beta+{\rm Re}(\mu)+\gamma(r)\eql0,$ where 
$\mu$ is given by $(1.7),$ $\omega=-a^*\alpha-\mu-1/2$ and $\zeta$ be chosen 
such that ${\rm Re}(\zeta)\eql\alpha$. Then if $\nu\notin{\cal E}_{\HHs},$ 
{\boldmath$H$}$(\LL_{\nu,r})$ can be represented in the form $(7.5),$ and if 
$\nu\in{\cal E}_{\HHs},$ then {\boldmath$H$}$(\LL_{\nu,r})$ is a subset of the 
right hand side of} (7.5).
\\\par
Finally we consider the case when $a^*>0,a^*_1<0$ and $a^*_2>0$. On the basis 
of Theorem 13, the following statement is proved similarly to that in the case 
$\delta=1$ in [20, Theorem 5.2].
\\\par
{\bf Theorem 14.} \ {\sl Let $a^*>0,$ $a^*_1<0,$ $a^*_2>0,$ 
$\alpha<1-\nu<\beta$ and let $1<r<\infty.$
\par
{\bf (a)} \ If $\nu\notin{\cal E}_{\HHs},$ or if $1<r\eql2,$ then 
{\boldmath$H$} is a one-to-one transform on $\LL_{\nu,r}$.
\par
{\bf (b)} \ Let
\begin{eqnarray}
&&\omega=a^*\eta-\Delta-\mu-\frac12,
\end{eqnarray}
where $\Delta$ and $\mu$ are given by $(1.6)$ and $(1.7)$ and let $\eta$ be 
chosen such that
\begin{eqnarray}
&&a^*{\rm Re}(\eta)\eqg\gamma(r)-2a^*_1\nu+\Delta+{\rm Re}(\mu),
\ {\rm Re}(\eta)>-\nu,
\end{eqnarray}
and let $\zeta$ be as
\begin{eqnarray}
&&{\rm Re}(\zeta)<\nu.
\end{eqnarray}
If $\nu\notin{\cal E}_{\HHs},$ then
\begin{eqnarray}
&&\mbox{\boldmath$H$}(\LL_{\nu,r})=\left(M_{-1/2-\omega/(2a^*_1)}
H_{2a^*_1,2a^*_1\zeta+\omega-1}L_{a^*,1/2-\eta+\omega/(2a^*_1)}\right)
\left(\LL_{1/2-{\rm Re}(\omega)/(2a^*_1)-\nu,r}\right),
\end{eqnarray}
where $M_\zeta$ is given by $(2.13)$. If $\nu\in{\cal E}_{\HHs},$ then 
{\boldmath$H$}$(\LL_{\nu,r})$ is a subset of the right hand side of} (7.15).
\\\par
{\bf Corollary 1.} \ {\sl Let $a^*>0,$ $a^*_1<0,$ $a^*_2>0,$ 
$-\infty\eql\alpha<1-\nu<\beta\eql\beta$ and let $1<r<\infty.$
\par
{\bf (a)} \ If $\nu\notin{\cal E}_{\HHs},$ or if $1<r\eql2,$ then 
{\boldmath$H$} is a one-to-one transform on $\LL_{\nu,r}$.
\par
{\bf (b)} \ Let $\omega$ be given by $(7.12)$ and let $\eta$ and $\zeta$ be 
chosen such that either of the following conditions holds
\begin{enumerate}
\item[{{\bf (i)}\ }]$a^*{\rm Re}(\eta)\eqg\gamma(r)-2a^*_1(1-\nu)+\Delta+
{\rm Re}(\mu),$ ${\rm Re}(\eta)\eqg\beta-1,$ ${\rm Re}(\zeta)\eql1-\beta,$ 
if $-\infty<\alpha<\beta<\infty;$
\item[{{\bf (ii)}\ }]$a^*{\rm Re}(\eta)\eqg\gamma(r)-2a^*_1\nu+
\Delta+{\rm Re}(\mu),$ ${\rm Re}(\eta)\eqg\beta-1,$ 
${\rm Re}(\zeta)\eql1-\beta,$ if $-\infty=\alpha<\beta<\infty;$
\item[{{\bf (iii)}\ }]$a^*{\rm Re}(\eta)\eqg\gamma(r)-2a^*_1(1-\nu)+
\Delta+{\rm Re}(\mu),$ ${\rm Re}(\eta)>-\nu,$ ${\rm Re}(\zeta)<\nu,$ if 
$-\infty<\alpha<\beta=\infty.$
\end{enumerate}
Then if $\nu\notin{\cal E}_{\HHs},$ {\boldmath$H$}$(\LL_{\nu,r})$ can be 
represented by the relation $(7.15),$ and if $\nu\in{\cal E}_{\HHs},$ 
{\boldmath$H$}$(\LL_{\nu,r})$ is a subset of the right hand side of} (7.15).
\\\par
{\bf Corollary 2.} \ {\sl Let $a^*>0,$ $a^*_1<0,$ $a^*_2>0$ and let 
$-\infty<\alpha<1-\nu<\beta<\infty,$ $1<r<\infty.$
\par
{\bf (a)} \ If $\nu\notin{\cal E}_{\HHs},$ or if $1<r\eql2,$ then 
{\boldmath$H$} is a one-to-one transform on $\LL_{\nu,r}$.
\par
{\bf (b)} \ Let $2a^*_1\alpha-a^*_2\beta+{\rm Re}(\mu)+\gamma(r)\eql0,$ 
where $\mu$ is given by $(1.7),$ and let $\zeta$ be chosen such that 
${\rm Re}(\zeta)\eql1-\beta$. Then if $\nu\notin{\cal E}_{\HHs},$ 
{\boldmath$H$}$(\LL_{\nu,r})$ can be represented by the relation $(7.15),$ 
and if $\nu\in{\cal E}_{\HHs},$ then {\boldmath$H$}$(\LL_{\nu,r})$ is a 
subset of the right hand side of} (7.15).
\\\\
\begin{center}
{\bf References}
\end{center}
\begin{enumerate}
\baselineskip=4.5mm
\item[{[1]\ }]V.M. Bhise and M. Dighe, On composition of integral operators 
with Fourier-type kernels, {\it Indian J. Pure Appl. Math.} {\bf 11}(1980), 
1183-1187.
\item[{[2]\ }]B.L.J. Braaksma, Asymptotic expansions and analytic continuation 
for a class of Barnes integrals, {\it Compositiones Math.} {\bf 15}(1964), 
239-341.
\item[{[3]\ }]Yu.A. Bruchkov, H.-J. Glaeske, A.P. Prudnikov and Vu Kim Tuan, 
{\it Multidimensional Integral Transformations}, Gordon and Breach, 
Philadelphia-Reading-Paris-Montreux-Tokyo-Melbourne, 1992.
\item[{[4]\ }]R.G. Buschman and H.M. Srivastava, Inversion formulas for the 
integral transformation with the $H$-function as kernel, {\it Indian J. 
Pure Appl. Math.} {\bf 6}(1975), 583-589.
\item[{[5]\ }]M. Dighe, Composition of fractional integral operator and an 
operator with Fourier-type kernel, {\it Bull Univ. Brasov} {\bf C20}(1978/79), 
3-8.
\item[{[6]\ }]A. Erd\'elyi, On some fractional transformations, {\it Rend. Sem. 
Mat. Univ. Politecn. Torino} {\bf 10}(1950), 217-234.
\item[{[7]\ }]A. Erd\'elyi, W. Magnus, F. Oberhettinger and F.G. Tricomi, 
{\it Higher Transcendental Functions}, Vol.\thinspace1, McGraw-Hill, 
New York-Toronto-London, 1953.
\item[{[8]\ }]C. Fox, The $G$ and $H$-functions as symmetrical Fourier 
kernels, {\it Trans. Amer. Math. Soc.} {\bf 98}(1961), 395-429.
\item[{[9]\ }]K.C. Gupta and P.K. Mittal, The $H$-function transform, 
{\it J. Austrial. Mth. Soc.} {\bf 11}(1970), 142-148.
\item[{[10]\ }]K.C. Gupta and P.K. Mittal, The $H$-function transform II, 
{\it J. Austrial. Mth. Soc.} {\bf 12}(1971), 445-450.
\item[{[11]\ }]S.L. Kalla, Integral operators involving Fox's $H$-function, 
{\it Acta Mexicana Ci. Tecn.} {\bf 3}(1969), 117-122.
\item[{[12]\ }]S.L. Kalla, Integral operators involving Fox's $H$-function II, 
{\it Notas Ci. Ser. M. Mat.} {\bf 7}(1969), 72-79.
\item[{[13]\ }]S.L. Kalla, Fractional integration operators involving 
generalized hypergeometric functions, {\it Univ. Nat. Tucum\'an Rev.}, 
Ser. A {\bf 20}(1970), 93-100.
\item[{[14]\ }]S.L. Kalla, On the solution of an integral equation involving a 
kernel of Mellin-Barnes type integral, {\it Kyungpook Math. J.} {\bf 12}(1972), 
93-101.
\item[{[15]\ }]S.L. Kalla, Operators of fractional integration, {\it Analytic 
Functions}, Kozubnik 1979, Lecture Notes in Math, {\bf 798}, 258-280, 
Springer, Belin, 1980.
\item[{[16]\ }]S.L. Kalla and V.S. Kiryakova, An $H$-function generalized 
fractional calculus based upon compositions of Erd\'elyi-Kober type operators 
in $L_p$, {\it Math. Japon.} {\bf 35}(1990), 1151-1171.
\item[{[17]\ }]R.N. Kesarwani, A pair of unsymmetrical Fourier kernels, 
{\it Trans. Amer. Math. Soc.} {\bf 115}(1965), 356-369.
\item[{[18]\ }]A.A. Kilbas, M. Saigo and S.A. Shlapakov, Integral transforms 
with Fox's $H$-function in spaces of summable functions, {\it Integral 
Transforms Spec. Funct.} {\bf 1}(1993), 87-103.
\item[{[19]\ }]A.A. Kilbas, M. Saigo and S.A. Shlapakov, Integral transforms 
with Fox's $H$-function in $L_{\nu,r}$-spaces, {\it Fukuoka Univ. Sci. Rep.} 
{\bf 23}(1993), 9-31.
\item[{[20]\ }]A.A. Kilbas, M. Saigo and S.A. Shlapakov, Integral transforms 
with Fox's $H$-function in $L_{\nu,r}$-spaces, II, {\it Fukuoka Univ. Sci. 
Rep.} {\bf 24}(1994), 13-38.
\item[{[21]\ }]A.A. Kilbas and S.A. Shlapakov, On a Bessel-type integral 
transformation and its compositions with integral and differential operators 
(Russian), {\it Dokl. Akad. Nauk Belarusi} {\bf 37}(1993), No.4, 10-14.
\item[{[22]\ }]A.A. Kilbas and S.A. Shlapakov, The composition of a 
Bessel-type integral operator with fractional integro-differentiation operators 
and its application to solving differential equations, {\it Differential 
Equations} {\bf 30}(1994), 235-246.
\item[{[23]\ }]A.A. Kilbas and S.A. Shlapakov, On an integral transform with 
Fox's $H$-function (Russian), {\it Dokl. Akad. Nauk Belarusi} {\bf 38}(1994), 
No.1, 12-15.
\item[{[24]\ }]V.S. Kiryakova, Fractional integration operators involving 
Fox's $H^{m,0}_{m,m}$-function, {\it C. R. Acad. Bulgare Sci.} {\bf 41}(1988), 
11-14.
\item[{[25]\ }]V.S. Kiryakova, {\it Generalized fractional Calculus and 
Applications}, Res. Notes Math., Vol.\thinspace301, Pitman, San Francisco 
{\it et alibi}, 1994.
\item[{[26]\ }]H. Kober, On fractional integrals and derivatives, {\it Quart. 
J. Math.}, {\it Oxford Ser.} {\bf 11}(1940), 193-211.
\item[{[27]\ }]E. Kr\"atzel, Integral transformations of Bessel-type, 
{\it General Functions and Operational Calculus}, (Proc. Conf., Varna, 1975), 
148-155, Bulg. Acad. Sci., Sofia, 1979.
\item[{[28]\ }]R.K. Kumbhat, On inversion formula for an integral transform, 
{\it Indian J. Pure Appl. Math.} {\bf 7}(1976), 368-375.
\item[{[29]\ }]A.M. Mathai and R.K. Saxena, {\it The $H$-Function with 
Applications in Statistics and Othe Disciplines}, Halsted Press [John Wiley 
\& Sons], New York-London-Sydney, 1978.
\item[{[30]\ }]A.C. McBride, Solution of hypergeometric integral equations 
involving generalized functions, {\it Proc. Edinburgh Math. Soc.} (2) 
{\bf 19}(1975), 265-285.
\item[{[31]\ }]A.C. McBride, {\it Fractional Calculus and Integral Transforms 
of Generalized Functions}, Res. Notes Math. {\bf 31}, Pitman, 
San Trancisco-London-Melbourne, 1979.
\item[{[32]\ }]A.C. McBride and W.J. Spratt, On the range and invertibility of 
a class of Mellin multiplier transforms, I, {\it J. Math. Anal. Appl.} 
{\bf 156}(1991), 568-587.
\item[{[33]\ }]A.C. McBride and W.J. Spratt, On the range and invertibility of 
a class of Mellin multiplier transforms, III, {\it Can. J. Math.} 
{\bf 43}(1991), 1323-1338.
\item[{[34]\ }]C. Nasim, An integral equation involving Fox's $H$-function, 
{Indian J. Pure Appl. Math.} {\bf 13}(1982), 1149-1162.
\item[{[35]\ }]T.R. Prabhaker, A class of integral equations with Gauss 
function in the kernels, {Math. Nachr.} {\bf 52}(1972), 71-83.
\item[{[36]\ }]A.P. Prudnikov, Yu.A. Brychkov and O.I. Marichev, 
{\it Integrals and Series}, Vol.\thinspace3, {\it More Special Functions}, 
Gordon \& Breach, New York {\it et alibi}, 1990.
\item[{[37]\ }]R.K. Raina and M. Saigo, Fractional calculus operators 
involving Fox's $H$-function in spaces $F_{p,\mu}$ and $F'_{p,\mu}$, 
{\it Recent Advances in Fractional Calculus}, 219-229, Global Publ., 
Souk Rapids (Minnesota), 1993.
\item[{[38]\ }]J. Rodriguez, J.J. Trujillo and M. Rivero, Operational 
fractional calculus of Kratzel integral transformation, {\it Differential 
equations} (Xanthi 1987), 613-620, Lect. Notes Pure Appl. Math., {\bf 118}, 
Dekker, New York, 1989.
\item[{[39]\ }]P.G. Rooney, A technique for studying the boundedness and 
extendability of certain types of operators, {\it Cand. J. Math.} 
{\bf 25}(1973), 1090-1102.
\item[{[40]\ }]P.G. Rooney, On integral transformations with $G$-function 
kernels, {\it Proc. Royal Soc. Edinburgh}, Sec. A {\bf 93}(1982/83), 265-297.
\item[{[41]\ }]M. Saigo, A remark on integral operators involving the Gauss 
hypergeometric functions, {\it Math. Rep. Kyushu Univ.} {\bf 11}(1978), 135-143.
\item[{[42]\ }]M. Saigo and H.-J. Glaeske, Fractional calculus operators 
involving the Gauss function in spaces $F_{p,\mu}$ and $F'_{p,\mu}$, 
{\it Math. Nachr.} {\bf 147}(1990), 285-306.
\item[{[43]\ }]M. Saigo, R.K. Raina and A.A. Kilbas, On generalized fractional 
calculus operators and their compositions with axisymmetric differential 
operator of the potential theory on spaces $F_{p,\mu}$ and $F'_{p,\mu}$, 
{\it Fukuoka Univ. Sci. Rep.} {\bf 23}(1993), 133-154.
\item[{[44]\ }]S.G. Samko, A.A. Kilbas and O.I. Marichev, {\it Fractional 
Integrals and Derivatives. Theory and Applications}, Gordon \& Breach, 
New York {\it et alibi}, 1993.
\item[{[45]\ }]R.K. Saxena, An inversion formula for a kernel involving a 
Mellin-Barnes type integral, {\it Proc. Amer. MAth. Soc.} {\bf 17}(1966), 
771-779.
\item[{[46]\ }]R.K. Saxena and R.K. Kumbhat, Integral operators involving 
$H$-function, {Indian J. Pure Appl. Math.} {\bf 5}(1974), 1-6.
\item[{[47]\ }]V.P. Saxena, Inversion formulae to certain integral equations 
involving $H$-functions, {Portugal Math.} {\bf 29}(1970), 31-42.
\item[{[48]\ }]S.A. Shlapakov, An integral transform involving Fox's 
$H$-function in the space of integrable functions (Russian), {\it Dokl. Akad. 
Nauk Belarus} {\bf 38}(1994), 14-18.
\item[{[49]\ }]R. Singh, An inversion formula for Fox $H$-transform, 
{\it Proc. Nat. Acad. Sci. India}, Sect A {\bf 40}(1970), 57-64.
\item[{[50]\ }]I.N. Sneddon, {\it Fractional Integration and Dual Integral 
Equations}, North Carolina State Coll., Appl. Math. Res. Group Raleigh, 
PSR-6, 1962.
\item[{[51]\ }]H.M. Srivastava and R.G. Buschman, Composition of fractional 
integral operators involving Fox's $H$-function, {\it Acta Mexicana Ci. Tecn.} 
{\bf 7}(1973), 21-28.
\item[{[52]\ }]H.M. Srivastava, K.C. Gupta and S.P. Goyal, {\it The 
$H$-functions of One and Two Variables with Applications}, South Asia Publ., 
New Delhi-Madras, 1982.
\item[{[53]\ }]Vu Kim Tuan, On the factorization of integral transformations 
of convolution type in the space $L_2^\Phi$ (Russian), {\it Dokl. Akad. Nauk 
Armenyan SSR} {\bf 83}(1986), 7-10.
\end{enumerate}
\end{document}